\documentclass{amsart}
\usepackage{amsmath,amsthm}
\usepackage{amsfonts,amssymb}

\usepackage[all]{xy}

\usepackage{enumerate}

\newtheorem{theorem}{Theorem}
\newtheorem{lemma}{Lemma}

\newtheorem{proposition}{Proposition}
\newtheorem{remark}{Remark}

\theoremstyle{remark}

\newcommand{\mR}{\mathbb{R}}
\newcommand{\mC}{\mathbb{C}}
\newcommand{\mN}{\mathbb{N}}
\newcommand{\mE}{\mathbb{E}}

\newcommand{\mS}{\mathbb{S}}

\newcommand{\cD}{\mathcal{D}}
\newcommand{\cM}{\mathcal{M}}
\newcommand{\cH}{\mathcal{H}}

\newcommand{\cF}{\mathcal{F}}
\newcommand{\cP}{\mathcal{P}}
\newcommand{\cC}{\mathcal{C}}

\newcommand{\cG}{\mathcal{G}}

\newcommand{\dD}{\textbf{D}}

\newcommand{\ux}{\underline{x}}

\newcommand{\uy}{\underline{y}}

\newcommand{\upx}{\partial_{\underline{x}}}

\begin{document}
 
\title[Dunkl operators and realizations of $\mathfrak{osp}(1|2)$]
{Dunkl operators and a family of realizations of $\mathfrak{osp}(1|2)$}
\author{H. De Bie}
\address{Department of Mathematical Analysis\\
Ghent University\\ Krijgslaan 281, 9000 Gent\\ Belgium.}
\email{Hendrik.DeBie@UGent.be}
\author{B. {\O}rsted}
\address{Department of Mathematical Sciences\\ University of Aarhus\\
Building 530, Ny Munkegade, DK 8000, Aarhus C\\ Denmark.}\email{orsted@imf.au.dk}

\author{P. Somberg}
\address{Mathematical Institute of Charles University\\
Sokolovsk\'a 83, 186 75 Praha\\ Czech Republic.}\email{somberg@karlin.mff.cuni.cz}

\author{V. Sou{\v{c}}ek}
\address{Mathematical Institute of Charles University\\
Sokolovsk\'a 83, 186 75 Praha\\ Czech Republic.}\email{soucek@karlin.mff.cuni.cz}

\date{\today}
\thanks{H. De Bie is a Postdoctoral Fellow of the Research Foundation - Flanders (FWO) and Courtesy Research Associate at the University of Oregon.} 

\begin{abstract}
In this paper, a family of radial deformations of the realization of the Lie superalgebra $\mathfrak{osp}(1|2)$ in the theory of Dunkl operators is obtained. This leads to a Dirac operator depending on 3 parameters. Several function theoretical aspects of this operator are studied, such as the associated measure, the related Laguerre polynomials and the related Fourier transform. For special values of the parameters, it is possible to construct the kernel of the Fourier transform explicitly, as well as the related intertwining operator.
\end{abstract}

\maketitle

\tableofcontents

\section{Introduction}
\label{preliminaries}

The notion of Howe duality (see \cite{MR0986027}) has received considerable attention lately. The basic framework where this duality is apparent is in the setting of classical harmonic analysis in $\mR^{m}$. Indeed, introducing the $O(m)$-invariant operators
\begin{eqnarray*}
\Delta &=& \sum_{i=1}^m \partial_{x_i}^2 \qquad \mbox{Laplace operator}\\
r^2 &=& \sum_{i=1}^m x_i^2\\
\mE &=&\sum_{i=1}^m x_i \partial_{x_i} \qquad \mbox{Euler operator},
\end{eqnarray*}
it is well-known that $\Delta$, $r^{2}$ and $\mE + \frac{m}{2}$ generate the Lie algebra $\mathfrak{sl}_{2}$, i.e.
\begin{eqnarray*}
\left[\Delta, r^2 \right] &=& 4 (\mE + \frac{m}{2})\\
\left[\Delta, \mE + \frac{m}{2} \right] &=& 2\Delta\\
\left[r^2, \mE + \frac{m}{2}\right] &=& -2 r^{2}.
\end{eqnarray*}
More importantly, the $\mathfrak{sl}_{2}$ action completely reduces the decomposition of the space $\cP$ of polynomials in $\mR^{m}$ to irreducible subspaces (the so-called spherical harmonics) under the action of the orthogonal group $O(m)$. In other words, the action of the dual pair $O(m) \times \mathfrak{sl}_{2}$ yields the complete multiplicity free decomposition of $\cP$. This decomposition is usually called Fischer decomposition and has, apart from its obvious importance in representation theory, also applications in e.g. the solution of the algebraic Dirichlet problem (see e.g. \cite{MR2031398} and references therein).

In the late eighties, Dunkl introduced in his seminal paper \cite{MR951883} the so-called Dunkl operators. These operators $T_{i}$, $i=1, \ldots, m$ are deformations of the classical partial derivatives by non-local operators which form a commutative family and allow to construct the Dunkl Laplacian $\Delta_{k} = \sum_{i=1}^m T_{i}^2$. This operator is no longer invariant under the whole orthogonal group but only under a finite reflection group $\cG$ (a finite subgroup of $O(m)$). Although this operator is much more complicated than the classical Laplace operator, Heckman noted in \cite{He} that the $\mathfrak{sl}_{2}$ relations are still valid in this case. Namely, one has
\begin{eqnarray*}
\left[\Delta_{k}, r^2 \right] &=& 4 (\mE + \frac{\mu}{2})\\
\left[\Delta_{k}, \mE + \frac{\mu}{2} \right] &=& 2\Delta_{k}\\
\left[r^2, \mE + \frac{\mu}{2}\right] &=& -2 r^{2}.
\end{eqnarray*}
where $\mu$ is a numerical parameter related to the group $\cG$ (for a precise definition we refer the reader to section \ref{DunklOps}).

The introduction of the Dunkl operator introduces a twist in the Fischer decomposition of the space of polynomials $\cP$ that preserves the radial structure but changes what happens on the sphere. That is, spherical harmonics for the classical Laplace operator are no longer necessarily spherical harmonics for the Dunkl Laplacian. For a discussion of the Fischer decomposition in both the classical and Dunkl case we refer the reader to \cite{Orsted}.

A natural question to ask now is whether one can also change the radial structure but nevertheless retain the $\mathfrak{sl}_{2}$ structure. This is taken up in  a recent preprint (see \cite{Orsted2}), where it was noted that the $\mathfrak{sl}_{2}$ relations also hold for the deformed operators $r^{2-a} \Delta_{k}$, $r^{a}$ and $\mE + \frac{a+m-2}{2}$, with $a$ in general a complex parameter, i.e.
\begin{eqnarray*}
\left[r^{2-a}\Delta_{k}, r^a \right] &=& 2 a (\mE + \frac{a +\mu -2}{2})\\
\left[r^{2-a}\Delta_{k}, \mE + \frac{a+\mu-2}{2} \right] &=& a r^{2-a}\Delta_{k}\\
\left[r^a, \mE + \frac{a +\mu-2}{2}\right] &=& -a r^{a}.
\end{eqnarray*}
In that paper, several other results pertaining to these deformed operators are obtained. In particular the associated Laguerre polynomials are developed for the study of the Fourier transform and the connection with the Gelfand-Gindikin program (for understanding a family of irreducible representations by using complex geometric methods) is elucidated. Two special cases are given a more detailed treatment, namely the classical case $a=2$ and the case $a=1$ (the study of which was introduced in \cite{MR2134314, MR2401813}).

Instead of considering the $\mathfrak{sl}_{2}$ Lie algebra generated by the three basic operators, one can also introduce two new operators, namely the Dirac operator and the vector variable. These operators enable us to factorize the Laplace operator and $r^{2}$. In the case of the orthogonal group, this yields the theory of Clifford analysis (see e.g. \cite{MR697564, MR1130821, MR1169463}). Indeed, by introducing the orthogonal Clifford algebra $\cC l_{0,m}$ one can refine $\mathfrak{sl}_{2}$ to the Lie superalgebra $\mathfrak{osp}(1|2)$. Considering an orthonormal basis $\{e_{i}\}$, $i=1, \ldots, m$ of $\mR^{m}$, we have that $\cC l_{0,m}$ is generated by this basis under the relations
\[
e_{i}e_{j} + e_{j} e_{i}= -2 \delta_{ij}.
\]
We then introduce the Dirac operator $\upx$ and the vector variable $\ux$ by
\begin{eqnarray*}
\upx &=& \sum_{i=1}^m e_i \partial_{x_i}, \qquad \mbox{Dirac operator}\\
\ux &=& \sum_{i=1}^m e_i x_i, \qquad \mbox{vector variable}.
\end{eqnarray*}
These operators satisfy
\begin{eqnarray*}
\upx^{2} &=& -\Delta\\
\ux^{2} &=& -r^{2}\\
\upx \ux + \ux \upx &=& -2 (\mE + \frac{m}{2})
\end{eqnarray*}
and generate a realization of $\mathfrak{osp}(1|2)$ (for the complete defining relations of this algebra, we refer the reader to theorem \ref{ospFamily}).
In this case, the classical Fischer decomposition given by the dual pair $O(m) \times \mathfrak{sl}_{2}$ can be refined using the dual pair $Spin(m) \times \mathfrak{osp}(1|2)$ if one lets the operators act on spinor-valued polynomials.

As the Dunkl operators are commutative, there also exists a Dunkl version of the Dirac operator, given by
\[
\cD_{k} = \sum_{i=1}^m e_i T_{i}
\]
which satisfies (see \cite{Orsted})
\[
\cD_{k} \ux + \ux \cD_{k} = -2(\mE + \frac{\mu}{2}).
\]

Let us now discuss the motivation and aims of the present paper. Firstly, finding a set of commuting operators $D_{i}$, $i=1,\ldots,m$ such that $\sum_{i=1}^{m} D_{i}^{2}= r^{2-a}\Delta_{k}$ would seriously facilitate the study of the harmonic analysis related to the operator $r^{2-a}\Delta_{k}$. For example, the Fourier transform related to $r^{2-a}\Delta_{k}$ is defined by the operator exponential (see \cite{Orsted2})
\[
\cF_{a} = e^{ \frac{i \pi (\mu+a-2)}{2a}} e^{\frac{i \pi}{2 a}(r^{2-a}\Delta_{k} - r^{a})}
\]
and the eigenfunctions and some of the properties are studied in detail there. However, except for the rank one case and the cases $a=1$ and $a=2$, the kernel of this integral transform is not explicitly known. In the non-Dunkl case the existence of the operators $D_{i}$ would yield a system of PDEs allowing to either determine the kernel (if one can explicitly solve the system) or at least obtain more information about its analytic behaviour (in a similar way as is done in e.g. \cite{deJ}).

Secondly, from the point of view of Clifford analysis, the existence of a Dirac operator related to $r^{2-a}\Delta_{k}$ would allow us to refine the $\mathfrak{sl}_{2}$ relations to $\mathfrak{osp}(1|2)$, hence enabling to extend the related function theory to the radially deformed case. Several interesting problems can then be addressed, such as determination of a fundamental solution, determination of the eigenfunctions and spectrum of the operator, existence of related Cauchy integral formulae etc.

As it turns out, it will not be possible to obtain a factorization of $r^{2-a}\Delta_{k}$ except for the special case $a=\pm 2$. We will however obtain an interesting new family of Dirac type operators, parametrized by three complex parameters $a,b$ and $c$ and defined as follows
\[
\dD = r^{1- \frac{a}{2}}\cD_{k} + b r^{-\frac{a}{2}-1} \ux + c r^{-\frac{a}{2}-1}\ux \mE.
\]
The square of this operator consists of $r^{2-a}\Delta_{k}$ plus some other terms. If one chooses the parameter $c=2/a -1$ the square of $\dD$ is scalar. For other values of $c$ the square of $\dD$ contains a bivector term. 

In other words, by further deforming the operator $r^{2-a}\Delta_{k}$ introduced in \cite{Orsted2}, one obtains an operator which can be factorized nicely. The properties of this operator are directly related to those of the classical Dirac operator via a generalized Kelvin transform. This Kelvin transform takes over the role played by the intertwining operator in the Dunkl case (see e.g. \cite{MR1273532}), although in our case we obtain an explicit formula for the operator.

The paper is organized as follows. In section \ref{DunklOps} we introduce the necessary background on Dunkl operators. In section \ref{SecFamily} we obtain the new family of Dirac operators and show that they generate the Lie superalgebra $\mathfrak{osp}(1|2)$. We explicitly calculate the square of $\dD$ and determine the measure related to $\dD$. In section \ref{SecLaguerre} we first discuss the analogues of polynomial solutions for $\dD$ and obtain a Fischer decomposition. Next we introduce the Laguerre polynomials and functions related to $\dD$. This allows us to introduce in section \ref{SecFourier} a new class of Fourier transforms, which are related to $\dD$. We determine the kernel explicitly when the square of $\dD$ is scalar and also discuss the very special case where $a=-2$. In section \ref{SecConclusions} we discuss several topics for further research.

\begin{remark}
The Dirac operator $\dD$ that we will introduce in section \ref{SecFamily} depends on three parameters $a, b$ and $c$. To simplify the reading, we will in the subsequent theorems not always mention the values that can be taken by these parameters. This can be verified easily by the reader when consulting the resulting formulae. As a general rule we can say that although in the basic definition all parameters can be complex, we will usually restrict ourselves to $a>0$ and $b$, $c$ real. Moreover, we always exclude the case $c=-1$.
Concerning the multiplicity function $k$ related to the Dunkl operators, we will restrict ourselves to $k > 0$ (although extensions are again possible). From section \ref{SecLaguerre} on, we will also exclude the parameter values given by formula (\ref{SingularLocus}). These values have to be compared with the singular values of $k$ described in \cite{MR1273532}.
\end{remark}

\section{Preliminaries on Dunkl operators}
\label{DunklOps}

Denote by $\langle .,. \rangle$ the standard Euclidean scalar product in $\mR^{m}$ and by $|x| = \langle x, x\rangle^{1/2}$ the associated norm. For $\alpha \in \mR^{m} - \{ 0\}$, the reflection $r_{\alpha}$ in the hyperplane orthogonal to $\alpha$ is given by
\[
r_{\alpha}(x) = x - 2 \frac{\langle \alpha, x\rangle}{|\alpha|^{2}}\alpha, \quad x \in \mR^{m}.
\]

A root system is a finite subset $R \subset \mR^{m}$ of non-zero vectors such that, for every $\alpha \in R$, the associated reflection $r_{\alpha}$ preserves $R$. We will assume that $R$ is reduced, i.e. $R \cap \mR \alpha = \{ \pm \alpha\}$ for all $\alpha \in R$. Each root system can be written as a disjoint union $R = R_{+} \cup (-R_{+})$, where $R_{+}$ and $-R_{+}$ are separated by a hyperplane through the origin. The subgroup $\cG \subset O(m)$ generated by the reflections $\{r_{\alpha} | \alpha \in R\}$ is called the finite reflection group associated with $R$. Following standard references (see e.g. \cite{deJ, MR1620515}) the roots in the root systems can be uniformly normalized in such a way that $\langle \alpha, \alpha\rangle = 2$ for all $\alpha\in R$. We assume such normalization throughout the article without further notice. For more information on finite reflection groups we refer the reader to \cite{Humph}.

A multiplicity function $k$ on the root system $R$ is a $\cG$-invariant function $k: R \rightarrow \mC$, i.e. $k(\alpha) = k(h \alpha)$ for all $h \in \cG$. We will denote $k(\alpha)$ by $k_{\alpha}$.

Fixing a positive subsystem $R_{+}$ of the root system $R$ and a multiplicity function $k$, we introduce the Dunkl operators $T_{i}$ associated to $R_{+}$ and $k$ by (see \cite{MR951883, MR1827871})
\[
T_{i} f(x)= \partial_{x_{i}} f(x) + \sum_{\alpha \in R_{+}} k_{\alpha} \alpha_{i} \frac{f(x) - f(r_{\alpha}(x))}{\langle \alpha, x\rangle}, \qquad f \in C^{1}(\mR^{m})
\]
where $\alpha_{i}$ is the $i$-th coordinate of $\alpha$.
An important property of the Dunkl operators is that they commute, i.e. $T_{i} T_{j} = T_{j} T_{i}$.

The Dunkl Laplacian is given by $\Delta_{k} = \sum_{i=1}^{m} T_i^2$, or more explicitly by
\[
\Delta_{k} f(x) = \Delta f(x) + 2 \sum_{\alpha \in R_{+}} k_{\alpha} \left( \frac{\langle \nabla f(x), \alpha \rangle}{\langle \alpha, x \rangle}  - \frac{f(x) - f(r_{\alpha}(x))}{\langle \alpha, x \rangle^{2}} \right)
\]
with $\Delta$ the classical Laplacian and $\nabla$ the gradient operator.

If we let $\Delta_{k}$ act on $|x|^2$ we find $\Delta_{k} |x|^2 = 2m + 4 \gamma = 2 \mu$, where $\gamma = \sum_{\alpha \in R_+} k_{\alpha}$. We call $\mu$ the Dunkl dimension, because most special functions related to $\Delta_{k}$ behave as if one would be working with the classical Laplace operator in a space with dimension $\mu$. We also denote by $\cH_\ell$ the space of Dunkl-harmonics of degree $\ell$, i.e. $\cH_\ell = \cP_\ell \cap \ker{\Delta_{k}}$ with $\cP_{\ell}$ the space of homogeneous polynomials of degree $\ell$. The space of Dunkl-harmonics of degree $\ell$ has the same dimension as the classical space of spherical harmonics of degree $\ell$.

It is possible to construct an intertwining operator $V_{k}$ connecting the classical derivatives $\partial_{x_{j}}$ with the Dunkl operators $T_{j}$ such that $T_{j} V_{k} = V_{k} \partial_{x_{j}}$ (see e.g. \cite{MR1273532}). Note that explicit formulae for $V_{k}$ are only known in a few special cases.

The operators $\Delta_{k}$, $r^{2}$ and $\mE + \frac{\mu}{2}$ satisfy the defining relations of the Lie algebra $\mathfrak{sl}_2$ (see \cite{He})
\begin{eqnarray*}
\left[\Delta_{k}, r^2 \right] &=& 4 (\mE + \frac{\mu}{2})\\
\left[\Delta_{k}, \mE + \frac{\mu}{2} \right] &=& 2\Delta_{k}\\
\left[r^2, \mE + \frac{\mu}{2}\right] &=& -2 r^{2}.
\end{eqnarray*}

We collect some basic properties of Dunkl operators in the following proposition.
\begin{proposition}
One has the following relations
\begin{equation*}
\begin{array}{ll}
(i)&T_{i}(f g) = (\partial_{i} f) g + f (T_{i}g), \qquad \mbox{if $f$ is invariant under $\cG$ and $f,g \in C^{1}(\mR^{m})$}\\
\vspace{-3mm}\\
(ii)&[x_{i}, \Delta_{k}] = - 2 T_{i}\\
\vspace{-3mm}\\
(iii)& \sum_{j=1}^{m} \left(x_{j} T_{j} + T_{j} x_{j} \right) = 2\mE + \mu\\
\vspace{-3mm}\\
(iv)& \partial_{r} T_{i} = T_{i} \partial_{r} + \dfrac{x_{i}}{r^{2}} \partial_{r} - \dfrac{1}{r} T_{i}.
\end{array} 
\end{equation*}
\label{BasicProps}
\end{proposition}

\begin{proof}
$(i)$ see \cite{MR951883}; $(ii)$ see \cite{MR951883}; $(iii)$ see \cite{He}.

Now we prove property $(iv)$.
As the Dunkl operators $T_{i}$ are operators of degree $-1$ and as the Euler operator can be written as $\mE = r \partial_{r}$, we have
\[
r\partial_{r} T_{i} = T_{i} r \partial_{r} - T_{i}.
\]
Applying $(i)$ then yields
\[
r \partial_{r} T_{i} = r T_{i} \partial_{r} + \frac{x_{i}}{r} \partial_{r} - T_{i}.
\]
Dividing both sides by $r$ gives the desired result.
\end{proof}

For the sequel, we also need to know that the operator $T_{i}x_{j} + x_{i}T_{j}$ is symmetric under $i \leftrightarrow j$. This is obtained in the following way
\begin{eqnarray*}
T_{i}x_{j} + x_{i}T_{j} &=& x_{i} \partial_{x_{j}} + x_{j} \partial_{x_{i}} + \sum_{\alpha \in R_{+}} k_{\alpha}\alpha_{j}x_{i} \frac{f(x) - f(r_{\alpha}(x))}{\langle \alpha, x \rangle}\\
&& + \sum_{\alpha \in R_{+}} k_{\alpha}\alpha_{i} \frac{x_{j}f(x) - [r_{\alpha}(x))]_{j}f(r_{\alpha}(x))}{\langle \alpha, x \rangle}\\
&=&  x_{i} \partial_{x_{j}} + x_{j} \partial_{x_{i}} + \sum_{\alpha \in R_{+}} k_{\alpha} (\alpha_{j}x_{i} + \alpha_{i}x_{j}) \frac{f(x)}{\langle \alpha, x \rangle}\\
&& - \sum_{\alpha \in R_{+}} k_{\alpha}(\alpha_{j} x_{i} + \alpha_{i} x_{j}-\langle\alpha,x\rangle \alpha_{i}\alpha_{j}) \frac{ f(r_{\alpha}(x))}{\langle \alpha, x \rangle},
\end{eqnarray*}
which is clearly symmetric in $i$ and $j$.

The weight function related to the root system $R$ and the multiplicity function $k$ is given by $w_k(x) = \prod_{\alpha \in R_{+}} |\langle \alpha, x\rangle |^{2 k_{\alpha}}$.
For suitably chosen functions $f$ and $g$ one then has the following property of integration by parts (see \cite{Du4})
\begin{equation}
\int_{\mR^{m}} (T_{i} f) g \;  w_{k}(x)dx = -\int_{\mR^{m}} f \left(T_{i}g\right) w_{k}(x)dx
\label{skewnessDunkl}
\end{equation}
with $dx$ the Lebesgue measure.

Recall that the classical Fourier transform is given by
\[
\cF_{\mbox{class}}(f)(y) = (2 \pi)^{-\frac{m}{2}} \int_{\mR^{m}} e^{-i\langle \ux,\uy\rangle} f(x) dx, \quad \langle \ux,\uy\rangle = \sum_{j=1}^{m}x_{j}y_{j}.
\]
There also exists a Fourier transform related to the set of Dunkl operators $T_{i}$ (see a.o. \cite{deJ}). This so-called Dunkl transform $\mathcal \cF_{k}: L^1(\mR^m, w_k(x)dx)\to C(\mR^m)$ is defined as follows
\[
 \mathcal \cF_{k} f(y):= c_{k}^{-1} \int_{\mR^m} f(x)\,D(x,-iy)\,w_k(x) dx \quad 
(y \in \mR^m)
\]
with $c_{k} = \int_{\mR^m} e^{-r^{2}/2}w_k(x) dx$ the Mehta constant related to $\cG$ and where $D(x,y)$ is the Dunkl kernel. This kernel is the unique solution of the system
\[
T_{i, x} D(x,y) = y_{i} D(x,y), \quad i=1, \ldots, m
\]
which is real-analytic in $\mR^{m}$ and satisfies $K(0,y)=1$.
The eigenfunctions of this transform are studied in a.o. \cite{Du4, MR1620515}. They are given by
\[
\phi_{j, \ell}^{k} = L_{j}^{\frac{\mu}{2} + \ell-1}(r^2) H_\ell \, e^{-r^{2}/2}, \quad H_{\ell} \in \cH_{\ell}
\]
and satisfy $\cF_{k}(\phi_{j, \ell}^{k}) = (-i)^{2j + \ell} \phi_{j, \ell}^{k}$.

 There also exists an exponential notation of both transforms
\begin{eqnarray*}
\cF_{\mbox{class}} &=&e^{ \frac{i \pi m}{4}} e^{\frac{i \pi}{4}(\Delta - r^{2})}\\
\cF_{k} &=& e^{ \frac{i \pi \mu}{4}} e^{\frac{i \pi}{4}(\Delta_{k} - r^{2})}
\end{eqnarray*}
In the Dunkl case, this is studied in depth in \cite{Said}. This notation is important for us, as it will provide us with the correct generalization of the Fourier transform to the new family of Dirac operators.

\section{A new family of realizations of $\mathfrak{osp}(1|2)$}
\label{SecFamily}

In this section we introduce a new class of Dirac operators that will generate $\mathfrak{osp}(1|2)$.
We try to stay as close as possible to the operators $r^{2-a} \Delta_{k}$, $r^{a}$ and $\mE + \frac{a+\mu-2}{2}$ introduced in \cite{Orsted2}.

We begin by factorizing the operator $r^{a}$. This is easily obtained by putting
\[
\ux_{a} = r^{\frac{a}{2}-1} \ux,
\]
because then $\ux_{a}^{2} = -r^{a}$.
Defining an $a$-deformed Dirac operator is less straight-forward. We take inspiration from the classical case $a=2$, where (see \cite{Orsted})
\[
\cD_{k} = -\frac{1}{2} [\ux, \Delta_{k}]
\]
and define
\[
\cD_{k,a} = -\frac{1}{2} [\ux_{a}, r^{2-a}\Delta_{k}].
\]
This commutator can be explicitly calculated as follows
\begin{eqnarray*}
\left[\ux_{a}, r^{2-a}\Delta_{k}\right] &=& r^{1-\frac{a}{2}} \ux \Delta_{k} - r^{2-a} \Delta_{k} r^{\frac{a}{2} -1} \ux\\
&=&  r^{1-\frac{a}{2}} \left[\ux, \Delta_{k}\right] - \left(\frac{a}{2} -1 \right)\left(\mu + \frac{a}{2} -3 \right) r^{2-a} r^{\frac{a}{2} -3}\ux\\
&& -2\left(\frac{a}{2} -1 \right)r^{2-a} r^{\frac{a}{2} -3} \mE \ux\\
&=& -2 r^{1-\frac{a}{2}} \cD_{k} - \left(\frac{a}{2} -1 \right)\left(\mu+\frac{a}{2} -3 \right) r^{-\frac{a}{2} -1}\ux\\
&& -2\left(\frac{a}{2} -1 \right)r^{-\frac{a}{2} -1} \ux \mE +2\left(\frac{a}{2} -1 \right)r^{-\frac{a}{2} -1} \ux\\
&=&-2 r^{1-\frac{a}{2}} \cD_{k} - \left(\frac{a}{2} -1 \right)\left(\mu+\frac{a}{2} -1 \right) r^{-\frac{a}{2} -1}\ux\\
&& -2\left(\frac{a}{2} -1 \right)r^{-\frac{a}{2} -1} \ux \mE
\end{eqnarray*}
where we have used proposition \ref{BasicProps}. Hence we obtain the following Ansatz for an $a$-deformed Dirac operator
\begin{equation}
\cD_{k,a} = r^{1- \frac{a}{2}}\cD_{k} + \frac{1}{2}\left(\frac{a}{2} -1 \right) \left(\frac{a}{2}+\mu-1 \right)r^{-\frac{a}{2}-1} \ux + \left( \frac{a}{2}-1\right)r^{-\frac{a}{2}-1}\ux \mE.
\label{aDirac}
\end{equation}
Note that this operator encompasses the classical Dirac operator ($a=2$, $k=0$) $\cD_{0,2} = \upx$ and the Dunkl Dirac operator ($a=2$) $\cD_{k,2} = \cD_{k}$. It turns out that we can work even slightly more general, by considering instead the operator
\begin{equation}
\dD = r^{1- \frac{a}{2}}\cD_{k} + b r^{-\frac{a}{2}-1} \ux + c r^{-\frac{a}{2}-1}\ux \mE.
\label{FamilyDirac}
\end{equation}
where $b$ and $c$ are now arbitrary complex numbers.

\vspace{3mm}
The operators $\dD$ and $\ux_{a}$ again generate a Lie superalgebra. This is the subject of the following theorem.
\begin{theorem}
The operators $\dD$ and $\ux_{a}$ generate a Lie superalgebra, isomorphic with $\mathfrak{osp}(1|2)$, with the following relations
\begin{equation}
\begin{array}{lll}
\{ \ux_{a}, \dD\} = -2(1+c) \left( \mE + \frac{\delta}{2}\right)&\quad&\left[\mE + \frac{\delta}{2}, \dD \right] = - \frac{a}{2} \dD\\
\vspace{-3mm}\\
\left[ \ux_{a}^{2}, \dD \right] = a(1+c)\ux_{a}&\quad&\left[\mE + \frac{\delta}{2}, \ux_{a} \right] =  \frac{a}{2} \ux_{a}\\
\vspace{-3mm}\\
\left[ \dD^{2}, \ux_{a} \right] = -a(1+c)\dD&\quad&\left[\mE + \frac{\delta}{2}, \dD^{2} \right] = - a \dD^{2}\\
\vspace{-3mm}\\
\left[ \dD^{2}, \ux_{a}^{2} \right] = 2a(1+c)^{2}\left( \mE + \frac{\delta}{2} \right)&\quad&\left[\mE + \frac{\delta}{2}, \ux_{a}^{2} \right] = a \ux_{a}^{2},
\end{array}
\end{equation}
where
\begin{equation}
\delta = \frac{a}{2} + \frac{2b + \mu -1}{1+c}.
\end{equation}
\label{ospFamily}
\end{theorem}

\begin{proof}
We start with the first anti-commutator
\begin{eqnarray*}
\{ \ux_{a}, \dD\} &=& \{ \ux_{a}, r^{1- \frac{a}{2}}\cD_{k} + b r^{-\frac{a}{2}-1} \ux + c r^{-\frac{a}{2}-1}\ux \mE \}\\
&=&\{ \ux_{a}, r^{1- \frac{a}{2}}\cD_{k}\} + c\{\ux_{a}, r^{-\frac{a}{2}-1}\ux \mE \} - 2b.
\end{eqnarray*}
The first term is calculated as
\[
\{ \ux_{a}, r^{1- \frac{a}{2}}\cD_{k}\}= \ux\cD_{k} +  r^{1- \frac{a}{2}}\cD_{k}  r^{ \frac{a}{2}-1}\ux = \{ \ux, \cD_{k}\} -\left( \frac{a}{2}-1 \right)
\]
and the second as
\[
\{\ux_{a},r^{-\frac{a}{2}-1}\ux \mE \} = -\mE + r^{-\frac{a}{2}-1}\ux\mE r^{\frac{a}{2}-1}\ux = -2\mE-\frac{a}{2}.
\]
Collecting everything then yields $\{ \ux_{a}, \dD\} = -2(1+c) \left( \mE + \frac{\delta}{2}\right)$.

Note that all relations in the second column are trivial. Using them and the previous calculation, we then obtain
\begin{eqnarray*}
\left[ \ux_{a}^{2}, \dD \right] &=& \ux_{a}^{2} \dD - \dD \ux_{a}^{2}\\
&=& \ux_{a}^{2} \dD + \ux_{a}\dD \ux_{a} - 2(1+c) \left( \mE + \frac{\delta}{2}\right) \ux_{a}\\
&=& -2(1+c) [\ux_{a},\mE]\\
&=& a(1+c) \ux_{a}.
\end{eqnarray*}
Similarly, we have
\begin{eqnarray*}
\left[ \dD^{2}, \ux_{a} \right] &=& \dD^{2} \ux_{a}-\ux_{a}\dD^{2}\\
&=& -\dD \ux_{a}\dD -\ux_{a} \dD^{2} -2(1+c)\dD\left( \mE + \frac{\delta}{2}\right)\\
&=& -2(1+c)[\dD,\mE]\\
&=& -a(1+c) \dD
\end{eqnarray*}
and finally
\begin{eqnarray*}
\left[ \dD^{2}, \ux_{a}^{2} \right] &=& \dD^{2} \ux_{a}^{2}-\ux_{a}^{2}\dD^{2}\\
&=& \dD \left(\ux_{a}^{2} \dD -a(1+c) \ux_{a} \right) -\ux_{a}^{2}\dD^{2}\\
&=&-a(1+c)\dD \ux_{a} + \left(\ux_{a}^{2} \dD -a(1+c) \ux_{a} \right) \dD - \ux_{a}^{2} \dD^{2}\\
&=&-a(1+c)\{ \ux_{a}, \dD\} \\
&=&2a(1+c)^{2}\left( \mE + \frac{\delta}{2} \right).
\end{eqnarray*}

Now, taking $\widetilde{a}$ and $\widetilde{1+c}$ such that $\widetilde{a}^{2}=a$ and $\widetilde{1+c}^{2} =1+c$, we can rescale the operators as follows
\[
\ux_{a} \rightarrow \ux_{a}/(\widetilde{a}\widetilde{1+c}), \quad \dD \rightarrow \dD/(\widetilde{a}\widetilde{1+c}), \quad \mE + \frac{\delta}{2} \rightarrow  (\mE + \frac{\delta}{2})/a
\]
which makes the isomorphism with $\mathfrak{osp}(1|2)$ explicit (see e.g. \cite{MR1773773}).
\end{proof}

Note that, although we have obtained a family of realizations of $\mathfrak{osp}(1|2)$, we have by no means shown that $\dD^{2} = -r^{2-a}\Delta_{k}$. Moreover, except for two special cases $a = \pm 2$ this will never be the case. In the following subsection we give the general result for the square of $\dD$, which is a complicated formula.

\begin{remark}
A special case of the operator $\dD$ has already been studied in \cite{MR2269930}. There the case $a=2, b=0$ and $c = -(1+ \alpha)$ is studied (for $k=0$) and the authors determine the eigenfunctions of this operator.
\end{remark}

\subsection{The square of $\dD$}

The deformed Dirac operator in formula (\ref{FamilyDirac}) is a vector-valued differential operator, i.e. $\dD$ is of the following form
\[
\dD= \sum_{i=1}^{m} D_{i}e_{i}.
\]
Introducing $l = 1-a/2$ to simplify notation, we are thus lead to consider the operators
\[
D_{i}= r^{l} T_{i} + b r^{l-2} x_{i} + c r^{l-1} x_{i} \partial_{r}, \qquad b,c, l \in \mC
\]
for $i = 1, \ldots, m$ where the $T_{i}$ are the Dunkl operators. The square of $\dD$ is then given by
\[
\dD^{2} = -\sum_{i=1}^{m} D_{i}^{2} + \sum_{i<j} e_{i}e_{j} (D_{i}D_{j} - D_{j} D_{i}),
\]
consisting of a scalar term $-\sum_{i=1}^{m} D_{i}^{2}$ and a bivector term $\sum_{i<j} e_{i}e_{j} (D_{i}D_{j} - D_{j} D_{i})$.

We calculate these two terms separately. First we simplify $D_{i} D_{j}$ by moving all Dunkl operators and derivatives to the right, and collecting all terms, symmetric under $i \leftrightarrow j$ in $\textbf{S}_{ij}$. We obtain
\begin{eqnarray*}
D_{i}D_{j} &=& \left(r^{l} T_{i} + b r^{l-2} x_{i} + c r^{l-1} x_{i} \partial_{r}\right)\left(r^{l} T_{j} + b r^{l-2} x_{j} + c r^{l-1} x_{j} \partial_{r}\right)\\
&=& r^{l} T_{i} r^{l} T_{j} +br^{l} T_{i}r^{l-2} x_{j} + cr^{l} T_{i}r^{l-1} x_{j} \partial_{r} + b r^{2l-2} x_{i}T_{j}
+b r^{2l-4} x_{i} x_{j} +\\
&& bc r^{2l-3} x_{i} x_{j} \partial_{r} + c r^{l-1} x_{i} \partial_{r}r^{l} T_{j} + b c r^{l-1} x_{i} \partial_{r}r^{l-2} x_{j} + c^{2} r^{l-1} x_{i} \partial_{r}r^{l-1} x_{j} \partial_{r}\\
&=& \textbf{S}_{ij} + l r^{2l-2} x_{i} T_{j} + r^{2l} T_{i} T_{j} + b(l-2) r^{2l-4}x_{i}x_{j} + b r^{2l-2}T_{i}x_{j} \\
&&+c(l-1)r^{2l-3}x_{i}x_{j}\partial_{r}+ c r^{2l-1} T_{i}x_{j} \partial_{r} + b r^{2l-2} x_{i}T_{j} + c l r^{2l-2} x_{i} T_{j}\\
&&+ c r^{2l-1}x_{i} \partial_{r}T_{j} + bc (l-2) r^{2l-4}x_{i}x_{j} 
 + bc r^{2l-3} x_{i}\partial_{r} x_{j}\\&& + c^{2}(l-1) r^{2l-3}x_{i}x_{j}\partial_{r} + c^{2}r^{2l-2}x_{i}\partial_{r}x_{j}\partial_{r}\\
&=& \textbf{S}_{ij} + (l + b+cl) r^{2l-2} x_{i}T_{j} + b r^{2l-2}T_{i}x_{j}+ c r^{2l-1} T_{i}x_{j}\partial_{r}\\&&+ c r^{2l-1} x_{i} \left(-T_{j} + T_{j}\partial_{r} + \frac{1}{r^{2}} x_{j}\partial_{r}\right)+bc r^{2l-3} x_{i} \left(\frac{x_{j}}{r} + x_{j}\partial_{r} \right)\\&& +c^{2} r^{2l-2} x_{i} \left(\frac{x_{j}}{r} + x_{j}\partial_{r} \right)\partial_{r}\\
&=& \textbf{S}_{ij}+ (l + b+cl-c) r^{2l-2} x_{i}T_{j} + b r^{2l-2}T_{i}x_{j} + c r^{2l-1} T_{i}x_{j}\partial_{r}\\&& +c r^{2l-1} x_{i}T_{j}\partial_{r}\\
&=&\textbf{S}_{ij} +(l + cl-c) r^{2l-2} x_{i}T_{j} + b r^{2l-2}\left(T_{i}x_{j} + x_{i}T_{j}\right)\\&&+ c r^{2l-1}\left(T_{i}x_{j} + x_{i}T_{j}\right)\partial_{r}\\
&=& \textbf{S}_{ij} + (l + cl - c) r^{2l-2} x_{i}T_{j},
\end{eqnarray*}
where we have used proposition \ref{BasicProps}, $(i)$ and $(iv)$, and also the fact that $T_{i}x_{j} + x_{i}T_{j}$ is symmetric.
Hence we conclude that $D_{i}D_{j} = D_{j} D_{i}$ if and only if
\[
l + cl - c =0
\]
or
\[
c = - \frac{l}{l-1} = \frac{2}{a}-1.
\]

We now calculate $\sum_{i=1}^{m} D_{i}^{2}$. We first have
\begin{eqnarray*}
D_{i}^{2} &=&  \left(r^{l} T_{i} + b r^{l-2} x_{i} + c r^{l-1} x_{i} \partial_{r}\right)  \left(r^{l} T_{i} + b r^{l-2} x_{i} + c r^{l-1} x_{i} \partial_{r}\right)\\
&=&r^{l}T_{i}r^{l}T_{i} + b^{2} r^{2l-4}x_{i}^{2} + c^{2} r^{l-1} x_{i} \partial_{r} r^{l-1} x_{i}\partial_{r} + b r^{l}T_{i} r^{l-2}x_{i}+b r^{2l-2}x_{i}T_{i}\\
&&+ c r^{l}T_{i} r^{l-1} x_{i}\partial_{r} + c r^{l-1} x_{i}\partial_{r} r^{l}T_{i} + bc r^{2l-3}x_{i}^{2} \partial_{r} + bc r^{l-1} x_{i}\partial_{r}r^{l-2}x_{i}\\
&=& b^{2}r^{2l-4} x_{i}^{2} + b r^{2l-2}x_{i}T_{i} + bc r^{2l-3}x_{i}^{2}\partial_{r} + r^{2l}T_{i}^{2} + l r^{2l-2}x_{i}T_{i}\\&& + c^{2} r^{2l-2}x_{i}\partial_{r}x_{i}\partial_{r}+ c^{2}(l-1) r^{2l-3}x_{i}^{2} \partial_{r} + b r^{2l-2}T_{i}x_{i} + b(l-2) r^{2l-4}x_{i}^{2}\\
&& + c(l-1) r^{2l-3} x_{i}^{2}\partial_{r} + c r^{2l-1}T_{i}x_{i}\partial_{r}+ cl r^{2l-2}x_{i}T_{i} + c r^{2l-1}x_{i}\partial_{r}T_{i} \\
&&+ bc(l-2)r^{2l-4} x_{i}^{2} + bc r^{2l-3} x_{i}\partial_{r}x_{i} \\
\end{eqnarray*}
\begin{eqnarray*}
&=&r^{2l}T_{i}^{2}+ \left(b^{2}+ b(l-2) + bc(l-2) \right) r^{2l-4}x_{i}^{2} + (b+l+ cl) r^{2l-2}x_{i}T_{i}\\
&& + \left(bc + c^{2}(l-1) + c(l-1) \right) r^{2l-3}x_{i}^{2}\partial_{r}+ c^{2}r^{2l-2}x_{i}\left(\frac{x_{i}}{r} + x_{i}\partial_{r} \right) \partial_{r}\\&& + b r^{2l-2}T_{i}x_{i}+ c r^{2l-1}T_{i}x_{i} \partial_{r}
+ c r^{2l-1}x_{i}\left(-\frac{1}{r} T_{i} + \frac{x_{i}}{r^{2}}\partial_{r} + T_{i}\partial_{r} \right)\\&&+bc r^{2l-3}x_{i}\left(\frac{x_{i}}{r} + x_{i}\partial_{r} \right) \\
&=&r^{2l}T_{i}^{2}+ \left(b^{2}+ b(l-2) + bc(l-2) + bc \right) r^{2l-4}x_{i}^{2} +(b+l+cl-c) r^{2l-2}x_{i}T_{i}+\\
&& + \left(2bc + c^{2}(l-1) + c(l-1) + c^{2} + c\right) r^{2l-3}x_{i}^{2}\partial_{r}+ c^{2}r^{2l-2}x_{i}^{2} \partial_{r}^{2} + b r^{2l-2}T_{i}x_{i}\\
&&+ c r^{2l-1}\left( T_{i}x_{i} + x_{i}T_{i}\right)\partial_{r}.
\end{eqnarray*}
Summing over $i = 1, \ldots, m$ then yields
\begin{eqnarray*}
\sum_{i=1}^{m} D_{i}^{2} &=& r^{2l} \Delta_{k} + \left(b^{2}+ b(l-2) + bc(l-1)  \right) r^{2l-2} + \left(2bc + c^{2}l + cl \right) r^{2l-1}\partial_{r}\\&& +  c^{2}r^{2l} \partial_{r}^{2}
+\sum_{i}(l+cl-c) r^{2l-2}x_{i}T_{i} + \sum_{i}br^{2l-2}\left(x_{i}T_{i} + T_{i}x_{i}\right)\\&& + \sum_{i}c r^{2l-1}\left( T_{i}x_{i} + x_{i}T_{i}\right)\partial_{r}\\
&=&r^{2l} \Delta_{k} + \left(b^{2}+ b(l-2) + bc(l-1) + b \mu \right) r^{2l-2}  +  \left(c^{2} + 2c\right)r^{2l} \partial_{r}^{2}\\&&+ \left(2bc + c^{2}l + cl + c\mu + 2b\right) r^{2l-1}\partial_{r}\\&&
+\sum_{i}(l+cl-c) r^{2l-2}x_{i}T_{i}.
\end{eqnarray*}
Hence $\sum_{i=1}^{n} D_{i}^{2} = r^{2-a} \Delta_{k}$ if and only if the following system of equations is satisfied
\begin{eqnarray*}
l &=& 1 - \frac{a}{2}\\
c (c+2) &=& 0\\
b\left(b + c (l-1) +(\mu+l-2)\right)&=&0\\
2b +  c (l+\mu) + lc^{2} + 2 bc &=& 0\\
l+cl-c&=&0.
\end{eqnarray*}
Note that the last equation is the same as the one guaranteeing the commutativity of the $D_{i}$. This system has exactly two solutions, namely
\begin{itemize}
\item $b=c=l=0$ and $a= 2$
\item $b=2-\mu$, $c=-2$, $l = 2$ and $a = -2$.
\end{itemize}

In the special case where the multiplicity function $k=0$ (so $\mu = m$), we obtain that $\sum_{i=1}^{n} D_{i}^{2} = r^{2-a} \Delta$ if and only if the following system of equations is satisfied
\begin{eqnarray*}
l &=& 1 - \frac{a}{2}\\
c (c+2) &=& 0\\
b\left(b + c (l-1) +(m+l-2)\right)&=&0\\
2b +l+  c (2l+m-1) + lc^{2} + 2 bc &=& 0,
\end{eqnarray*}
i.e. the last two equations in the previous system are merged. In this case, 
the system has exactly 4 solutions, namely
\begin{itemize}
\item $b=c=l=0$ and $a= 2$
\item $b=0$, $c=-2$, $l = 2-m$ and $a=2m-2$
\item $b=m-2$, $c=0$, $l= 4-2m$ and $a=2(2m-3)$
\item $b=2-m$, $c=-2$, $l = 2$ and $a = -2$
\end{itemize}
where now only the first and the last solution are compatible with the commutativity of the $D_{i}$. 

\vspace{3mm}
We summarize our results in the following theorem:
\begin{theorem}
The deformed Dunkl Dirac operator $\dD$ factorizes $r^{2-a}\Delta_{k}$ if and only if $a=2$ or $a=-2$. In the case $a=2$ this is the classical Dunkl Dirac operator $\cD_{k}$ satisfying $\cD_{k}^{2} = - \Delta_{k}$. In the case $a=-2$ this is the operator
\[
\cD_{k,-2} = r^{2}\cD_{k} - \left(\mu-2 \right) \ux -2 \ux \mE
\]
which satisfies $\cD_{k,a}^{2} = -r^{4}\Delta_{k}$.
\end{theorem}

The relation between $\cD_{k}$ and $\cD_{k,-2}$ will be discussed in section \ref{KelvinSection}.

For general values of $b$ and $c$, the square of $\dD$ is a complicated operator, given by
\begin{eqnarray*}
\dD^{2} &=&-r^{2-a} \Delta_{k} -b\left(b -1-\frac{a}{2}(1+c)+  \mu \right) r^{-a}\\
&& - \left(2bc + (c^{2} + c)(1-\frac{a}{2}) + c\mu + 2b\right) r^{1-a}\partial_{r}\\
&& -  \left(c^{2} + 2c\right)r^{2-a} \partial_{r}^{2}-(1-\frac{a}{2}(1+c)) r^{-a} \sum_{i} x_{i}T_{i}\\
&& + (1-\frac{a}{2}(1+c)) r^{-a}  \sum_{i<j} e_{i}e_{j} ( x_{i}T_{j} - x_{j}T_{i})
\end{eqnarray*}
and if $c= 2/a-1$ by
\begin{eqnarray*}
\dD^{2} &=&-r^{2-a} \Delta_{k} -b\left(b -1-\frac{a}{2}(1+c)+  \mu \right) r^{-a}\\
&& - \left(2bc + (c^{2} + c)(1-\frac{a}{2}) + c\mu + 2b\right) r^{1-a}\partial_{r} -  \left(c^{2} + 2c\right)r^{2-a} \partial_{r}^{2}.
\end{eqnarray*}

\subsection{The measure associated to $\dD$}

We want to associate a measure to the operator $\dD$ in such a way that we can perform integration by parts. Concretely, as the operator $\dD$ can be written as $\dD = \sum_{i=1}^{m} D_{i} e_{i}$ with
\[
D_{i} = r^{1-\frac{a}{2}} T_{i} + b r^{-1-\frac{a}{2}} x_{i} + c r^{-\frac{a}{2}} x_{i} \partial_{r}
\]
we want to determine a radial function $h(r)$ such that for all $i=1,\ldots,m$
\begin{equation}
\int_{\mR^{m}} (D_{i}f) g h(r) w_{k}(x)dx = -\int_{\mR^{m}} f \left(D_{i}g\right) h(r) w_{k}(x)dx
\label{partIntAnsatz}
\end{equation}
for all scalar functions $f$ and $g$ such that the integrals exist and formula (\ref{skewnessDunkl}) can be applied.

We start by calculating the left-hand side, yielding
\begin{eqnarray*}
\int_{\mR^{m}} (D_{i}f) g h(r) w_{k}(x)dx &=& \int_{\mR^{m}}  r^{1-\frac{a}{2}}(T_{i}f) g h(r) w_{k}(x)dx\\
&& + b \int_{\mR^{m}}  r^{-1-\frac{a}{2}} x_{i}f g h(r) w_{k}(x)dx\\
&&+ c\int_{\mR^{m}} r^{-\frac{a}{2}} x_{i}(\partial_{r}f) g h(r) w_{k}(x)dx.
\end{eqnarray*}
The first integral can be rewritten, using formula (\ref{skewnessDunkl}), as 
\begin{eqnarray*}
\int_{\mR^{m}}  r^{1-\frac{a}{2}}(T_{i}f) g h(r) w_{k}(x)dx &=& -\int_{\mR^{m}}  f T_{i} \left(r^{1-\frac{a}{2}} g h(r)\right) w_{k}(x) dx\\
&=& -\int_{\mR^{m}}  f (T_{i} g) r^{1-\frac{a}{2}} h(r) w_{k}(x) dx\\
&& - \int_{\mR^{m}}  f g \left(\partial_{r} h(r) \right)r^{-\frac{a}{2}} x_{i} w_{k}(x) dx\\
&&- (1-\frac{a}{2})\int_{\mR^{m}}  f g h(r) r^{-1-\frac{a}{2}}x_{i} w_{k}(x) dx.
\end{eqnarray*}
Now we consider the third integral. First note that
\begin{equation}
\int_{\mR^{m}} \left(\partial_{r}f\right) g dx = -\int_{\mR^{m}} f\left(\partial_{r}g\right)  dx + (1-m) \int_{\mR^{m}} f g r^{-1} dx.
\label{skewnessDr}
\end{equation}
Using this, we obtain
\begin{eqnarray*}
\int_{\mR^{m}} r^{-\frac{a}{2}} x_{i}(\partial_{r}f) g h(r) w_{k}(x)dx &=& -\int_{\mR^{m}} f \partial_{r} \left(r^{-\frac{a}{2}} x_{i} g h(r) w_{k}(x) \right)dx\\
&&+ (1-m)\int_{\mR^{m}} r^{-\frac{a}{2}-1} x_{i}f g h(r) w_{k}(x)dx\\
&=&-\int_{\mR^{m}} f \left(\partial_{r} g \right) r^{-\frac{a}{2}} x_{i} h(r) w_{k}(x) dx\\
&& - \int_{\mR^{m}} f g \left(\partial_{r}h\right) r^{-\frac{a}{2}} x_{i} w_{k}(x) dx\\
&&+ (\frac{a}{2}-\mu)\int_{\mR^{m}} f g r^{-\frac{a}{2}-1} x_{i} h(r) w_{k}(x)dx.
\end{eqnarray*}
Combining the three integrals and comparing with the right-hand side in equation (\ref{partIntAnsatz}), we obtain the following differential equation for $h(r)$
\[
(1+c) r \frac{d}{dr} h = \left(2b -1+ \frac{a}{2} - c( \mu - \frac{a}{2})\right) h.
\]
We exclude the case $c = -1$ as it leads to singular behaviour. For all other values of $a,b$ and $c$, the solution is, up to a constant, given by
\[
h(r) = r^{\frac{a}{2} + \frac{2b-1-c\mu}{1+c}}.
\]
Note that the entire measure associated to $\dD$ is now given by $h(r)  w_{k}(x)dx$ and that its radial part is given by
\[
r^{\delta -1} dr
\]
which indicates that the parameter $\delta$ (appearing in theorem \ref{ospFamily}) behaves as the dimension describing the system.

Summarizing, we have obtained the following proposition.
\begin{proposition}
If $c \neq -1$, then for suitable differentiable functions $f$ and $g$ the following holds for all $i=1, \ldots, m$
\[
\int_{\mR^{m}} (D_{i}f) g h(r) w_{k}(x)dx = -\int_{\mR^{m}} f \left(D_{i}g\right) h(r) w_{k}(x)dx
\]
with $h(r) = r^{\frac{a}{2} + \frac{2b-1-c\mu}{1+c}}$, provided the integrals exist.
\label{partIntProp}
\end{proposition}

\section{Special functions related to $\dD$}
\label{SecLaguerre}

\subsection{The Fischer decomposition and the dual pair related to $\dD$}

We first consider null-solutions of $\dD$. Recall that polynomial null-solutions of the Dunkl-Dirac operator $\cD_{k}$ are the so-called Dunkl-monogenics (see \cite{Orsted}). More precisely, the space $\cM_{\ell}$ of Dunkl-monogenics of degree $\ell$ is the space of homogeneous polynomials of degree $\ell$ which are in the kernel of $\cD_{k}$, i.e. $\cM_{\ell} = \ker{\cD_{k}}\cap ( \cP_{\ell}\otimes \mS)$. In this notation $\mS$ is a representation of the Clifford algebra $\cC l_{0,m}$. Possible choices include the (irreducible) spinor spaces or the whole Clifford algebra itself. Note that clearly $\cM_{\ell} \subset \cH_{\ell}\otimes \mS$.

We now want to find the analogues of these solutions for the new operator $\dD$. Due to the form of $\dD$, it makes sense to propose solutions of the form
\[
f = r^{\beta} M_{\ell}, \qquad M_{\ell} \in \cM_{\ell}
\]
with $\beta$ to be determined. Expressing $\dD f = 0$ we obtain
\[
(\beta + b + c(\beta+\ell)) r^{\beta-2}\ux M_{\ell} =0
\]
yielding
\[
\beta = \beta_{\ell} = - \frac{b + c\ell}{1+c}.
\]
It is important to note that $\beta$ depends on the degree of the Dunkl-monogenic considered. Hence the space $r^{\beta_{\ell}} \cM_{\ell}$ takes over the role played by the space of Dunkl monogenics $\cM_{\ell}$.

Next we need the following basic lemma.
\begin{lemma}
One has the following relations
\begin{eqnarray*}
\dD\left( \ux_{a}^{2t} r^{\beta_{\ell}} M_{\ell} \right) &=& -  (1+c) a t \, \ux_{a}^{2t-1} r^{\beta_{\ell}}M_{\ell}\\
\dD \left(\ux_{a}^{2t+1} r^{\beta_{\ell}} M_{\ell} \right)&=& - (1+c) (\gamma_{\ell} + at) \, \ux_{a}^{2t} r^{\beta_{\ell}}M_{\ell}
\end{eqnarray*}
with $M_{\ell} \in \cM_{\ell}$, $\gamma_{\ell}  = 2 \beta_{\ell} + 2\ell + \delta$ and $t \in \mN$.
\label{FischerLemma}
\end{lemma}

\begin{proof}
We have, using theorem \ref{ospFamily},
\begin{eqnarray*}
\dD \left(\ux_{a} r^{\beta_{\ell}} M_{\ell} \right) &=& - \ux_{a}\dD \left( r^{\beta_{\ell}} M_{\ell} \right) -2(1+c)\left( \mE + \frac{\delta}{2}\right) r^{\beta_{\ell}} M_{\ell}\\
&=& -(1+c)(2\beta_{\ell} + 2\ell + \delta)r^{\beta_{\ell}} M_{\ell}. 
\end{eqnarray*}
The general formulae follow easily using induction.
\end{proof}

We define the space of homogeneically shifted Dunkl monogenics of degree $\ell$ by
\[
\cM_{\ell}^{\beta_{\ell}} = r^{\beta_{\ell}}\cM_\ell.
\]
So clearly $\dD \cM_{\ell}^{\beta_{\ell}} =0$. Starting from this space, we can now generate an infinite-dimensional representation of the Lie superalgebra $\mathfrak{osp}(1|2)$ as follows
\[
\xymatrix{\cM_{\ell}^{\beta_{\ell}} \ar@<1ex>[r]^{\ux_{a}} \ar@(dl,dr)_{\mE + \frac{\delta}{2}}   & \ux_{a} \cM_{\ell}^{\beta_{\ell}} \ar@(dl,dr)_{\mE + \frac{\delta}{2}} \ar@<1ex>[r]^{\ux_{a}} \ar@<1ex>[l]^{\dD}   & \ux_{a}^2 \cM_{\ell}^{\beta_{\ell}} \ar@(dl,dr)_{\mE + \frac{\delta}{2}} \ar@<1ex>[r]^{\ux_{a}} \ar@<1ex>[l]^{\dD} &\ux_{a}^3 \cM_{\ell}^{\beta_{\ell}} \ar@(dl,dr)_{\mE + \frac{\delta}{2}} \ar@<1ex>[r]^{\ux_{a}} \ar@<1ex>[l]^{\dD} &\ux_{a}^4 \cM_{\ell}^{\beta_{\ell}}  \ar@(dl,dr)_{\mE + \frac{\delta}{2}} \ar@<1ex>[r]^{\ux_{a}} \ar@<1ex>[l]^{\dD} &\ar@<1ex>[l]^{\dD} \ldots
}
\]
provided that $\gamma_{\ell}/a \not \in -\mN$ (see lemma \ref{FischerLemma}). This means that
\begin{equation}
\frac{\gamma_{\ell}}{a} = \frac{1}{2} + \frac{ \mu -1 + 2l}{a(1+c)} \not \in -\mN.
\label{SingularLocus}
\end{equation}
Note that this condition is e.g. always fulfilled if $a(1+c) >0$ and the multiplicity function $k >0$.

In the classical case $a=2$, $b=c=0$ one can put all these infinite-dimensional representations in one scheme as follows
\[
\xymatrix@=11pt{\cP_0 \otimes \mS  &  \cP_1 \otimes \mS & \cP_2 \otimes \mS & \cP_3 \otimes \mS \ar@{=}[d] & \cP_4 \otimes \mS & \cP_5 \otimes \mS  &\ldots \\
\cM_0 \ar[r] & \ux \cM_0 \ar[r] & \ux^2 \cM_0 \ar[r] & \ux^3 \cM_0 \ar @{} [d] |{\oplus}
 \ar[r] & \ux^4 \cM_0 \ar[r] & \ux^5 \cM_0 &\ldots\\
&\cM_1 \ar[r] & \ux \cM_1 \ar[r] & \ux^2 \cM_1 \ar @{} [d] |{\oplus}
 \ar[r] & \ux^3 \cM_1 \ar[r] & \ux^4 \cM_1  &\ldots\\
&&\cM_2 \ar[r] & \ux \cM_2 \ar @{} [d] |{\oplus}
 \ar[r] & \ux^2 \cM_2 \ar[r] & \ux^3 \cM_2  &\ldots\\
&&&\cM_3 \ar[r] & \ux \cM_3 \ar[r] & \ux^2 \cM_3  &\ldots\\
&&&&\cM_4 \ar[r] & \ux \cM_4   &\ldots\\
&&&&&\cM_5&\ldots
}
\]
Each column now yields the decomposition of the space of $\mS$-valued homogeneous polynomials of a certain degree into Dunkl monogenics.
If $\mS$ is an irreducible representation of $\cC l_{0,m}$, then we have obtained the Fischer decomposition of $\cP \otimes \mS$ under the action of the dual pair $Spin(m) \times \mathfrak{osp}(1|2)$ (when the multiplicity function $k=0$).

The fact that one can collect particular weight spaces generated by several homogeneically shifted monogenic functions according to the same value of $\mE+\delta/2$ is not generically true in the case of the deformed operator $\dD$. This is only the case when $\beta_{\ell-1} + (\ell-1) + a/2 = \beta_{\ell} +\ell$ or equivalently $c = 2/a-1$ (note that this is also the special case when the square of $\dD$ is scalar). In this case, the scheme is given by
\[
\xymatrix@=11pt{\cP_{\beta_{0}} \otimes \mS  &  \cP_{\beta_{1}+1} \otimes \mS & \cP_{\beta_{2}+2} \otimes \mS & \cP_{\beta_{3}+3} \otimes \mS \ar@{=}[d] & \cP_{\beta_{4}+4} \otimes \mS & \cP_{\beta_{5}+5} \otimes \mS  &\ldots \\
\cM_{0}^{\beta_{0}} \ar[r] & \ux_a \cM_{0}^{\beta_{0}} \ar[r] & \ux_a^2 \cM_{0}^{\beta_{0}} \ar[r] & \ux_a^3 \cM_{0}^{\beta_{0}} \ar @{} [d] |{\oplus}
 \ar[r] & \ux_a^4 \cM_{0}^{\beta_{0}} \ar[r] & \ux_a^5 \cM_{0}^{\beta_{0}} &\ldots\\
&\cM_{1}^{\beta_{1}} \ar[r] & \ux_a \cM_{1}^{\beta_{1}} \ar[r] & \ux_a^2 \cM_{1}^{\beta_{1}} \ar @{} [d] |{\oplus}
 \ar[r] & \ux_a^3 \cM_{1}^{\beta_{1}} \ar[r] & \ux_a^4 \cM_{1}^{\beta_{1}}  &\ldots\\
&&\cM_{2}^{\beta_{2}} \ar[r] & \ux_a \cM_{2}^{\beta_{2}} \ar @{} [d] |{\oplus}
 \ar[r] & \ux_a^2\cM_{2}^{\beta_{2}}  \ar[r] & \ux_a^3\cM_{2}^{\beta_{2}}   &\ldots\\
&&&\cM_{3}^{\beta_{3}}  \ar[r] & \ux_a \cM_{3}^{\beta_{3}}  \ar[r] & \ux_a^2 \cM_{3}^{\beta_{3}}   &\ldots\\
&&&&\cM_{4}^{\beta_{4}}  \ar[r] & \ux_a \cM_{4}^{\beta_{4}}    &\ldots\\
&&&&&\cM_{5}^{\beta_{5}} &\ldots
}
\]
In this scheme, spaces of the type $\cP_{\beta_{k}+k} \otimes \mS$ are defined by
\[
\cP_{\beta_{k}+k} \otimes \mS = \bigoplus_{j=0}^{k} \ux_{a}^{k-j}\cM_{j}^{\beta_{j}}
\]
where all summands have the same homogeneity, as measured by the Euler operator.

\subsection{Laguerre polynomials related to $\dD$}

In this subsection and the rest of the paper, we will always consider functions (and Dunkl monogenics) taking values in the full Clifford algebra $\cC l_{0,m}$.

We start with the following technical result.
\begin{lemma}[gauging of $\dD$]
One has the following operator equality:
\[
e^{r^{a}/a} \dD e^{-r^{a}/a} = \dD - (1+c)\ux_{a}.
\]
\label{gaugeD}
\end{lemma}
\begin{proof}
We calculate this identity as follows
\begin{eqnarray*}
e^{r^{a}/a} \dD e^{-r^{a}/a} &=& e^{r^{a}/a} \left( r^{1- \frac{a}{2}}\cD_{k} + b r^{-\frac{a}{2}-1} \ux + c r^{-\frac{a}{2}-1}\ux \mE \right) e^{-r^{a}/a}\\
&=&e^{r^{a}/a} r^{1- \frac{a}{2}}  \left(\upx e^{-r^{a}/a} \right) + r^{1- \frac{a}{2}}\cD_{k} + b r^{-\frac{a}{2}-1} \ux\\
&&  + c e^{r^{a}/a} r^{-\frac{a}{2}-1}\ux \left( \mE  e^{-r^{a}/a} \right) + c r^{-\frac{a}{2}-1}\ux \mE\\
&=& \dD- (1+c)\ux_{a}.
\end{eqnarray*}
\end{proof}

We now have all the tools necessary to introduce Laguerre polynomials and functions related to the operator $\dD$. 

To each Dunkl monogenic $M_{\ell} \in \cM_{\ell}$ we associate two sets of functions, namely the set of functions $\{\psi_{t,\ell}\}$ and the set of functions $\{\phi_{t,\ell}\}$ as follows:
\begin{equation}
\psi_{t,\ell} = \left(\dD - 2(1+c) \ux_{a} \right)^{t} r^{\beta_{\ell}} M_{\ell}, \qquad t \in \mN 
\label{defPsi}
\end{equation}
and
\begin{equation}
\phi_{t,\ell} = \left(\left(\dD - 2(1+c) \ux_{a} \right)^{t} r^{\beta_{\ell}} M_{\ell}\right) e^{-r^{a}/a} = \left(\dD - (1+c) \ux_{a} \right)^{t} r^{\beta_{\ell}} M_{\ell} e^{-r^{a}/a}, \, t \in \mN
\end{equation}
where the last equality follows from lemma \ref{gaugeD}.

The set of functions $\psi_{t,\ell}$ is a generalization of the so-called Clifford-Hermite polynomials introduced by Sommen in \cite{MR926831}. The set $\phi_{t,\ell} = \psi_{t,\ell}e^{-r^{a}/a}$ is then a generalization of the Clifford-Hermite functions.
These functions are very important, as they will turn out to be eigenfunctions of a generalized Fourier transform.
Now we proceed by obtaining the basic properties of these new functions.

\begin{theorem}[Differential equation]
\label{DiffequationCHpolys}
$\psi_{t,\ell}$ is a solution of the following differential equation:
\[
\dD^{2} \psi_{t,\ell} -2(1+c)  \ux_{a} \dD \psi_{t,\ell} - C(t,\ell) \psi_{t,\ell} = 0
\]
with 
\[
 \left\{ \begin{array}{l} 
C(2t,\ell) =2(1+c)^{2}at\\
C(2t+1,\ell) =2(1+c)^{2}(\gamma_{\ell} + at).
\end{array}
\right.
\]
\end{theorem}

\begin{proof}
We will in fact prove something slightly stronger, namely
\begin{equation}
\begin{array}{lll}
\dD \psi_{2t,\ell} &=&2(1+c)^{2}at \, \psi_{2t-1,\ell}\\
\dD \psi_{2t+1,\ell} &=&2(1+c)^{2}(\gamma_{\ell} + at) \, \psi_{2t,\ell}.
\end{array}
\label{pdehermite}
\end{equation}
The theorem then immediately follows by acting on (\ref{pdehermite}) with $\dD  -2(1+c)  \ux_{a}$.

Using formula (\ref{defPsi}) and lemma \ref{FischerLemma}, it is easy to see that we can expand the functions $\psi_{t,\ell}$ as follows
\[
\begin{array}{lll}
\psi_{2t,\ell} &=& \sum_{i=0}^t b_{2i}^{2t} \ux_{a}^{2i} r^{\beta_{\ell}}M_\ell\\
\psi_{2t+1,\ell} &=& \sum_{i=0}^t b_{2i+1}^{2t+1} \ux_{a}^{2i+1} r^{\beta_{\ell}} M_\ell.
\end{array}
\]
The recursion relation $\psi_{t+1,\ell} =\left(\dD - 2(1+c) \ux_{a} \right)\psi_{t,\ell}$ then yields the following relation among the coefficients
\[
\begin{array}{lll}
b_{2i}^{2t} &=& -(1+c)(\gamma_{\ell} + a i) b_{2i+1}^{2t-1} -2(1+c) b_{2i-1}^{2t-1} \\
b_{2i+1}^{2t+1} &=& -a (1+c)(i+1)b_{2i+2}^{2t} -2(1+c)b_{2i}^{2t}.
\end{array}
\]

In terms of the coefficients $b_j^i$, formula (\ref{pdehermite}) takes the following form
\[
\begin{array}{lll}
 i b_{2i}^{2t} &=& -2(1+c) t b_{2i-1}^{2t-1}\\
(\gamma_{\ell} + a i ) b_{2i+1}^{2t+1} &=& -2(1+c) (\gamma_{\ell} + a t)b_{2i}^{2t}.
\end{array}
\]
which can now be proven using induction. Indeed, it is easy to check the theorem for $t=0, 1$. So suppose that formula (\ref{pdehermite}) holds for $\psi_{t,k}$, $t \leq 2s$. We show that it also holds for $t=2s+1$. We have
\begin{eqnarray*}
(\gamma_{\ell} + ai) b_{2i+1}^{2s+1} &=& -(\gamma_{\ell} + a i) (1+c)( 2b_{2i}^{2s} + a(i+1)b_{2i+2}^{2s}  )\\
&=&-(\gamma_{\ell} + a i) (1+c)(  2b_{2i}^{2s}-2a (1+c)s b_{2i+1}^{2s-1} )\\
&=& -2(\gamma_{\ell} + a i) (1+c)b_{2i}^{2s} - 2a (1+c)s (b_{2i}^{2s}+2 (1+c) b_{2i-1}^{2s-1})\\
&=& -2(\gamma_{\ell} + a t) (1+c)b_{2i}^{2s} - 2a (1+c)i b_{2i}^{2s} -4a (1+c)^{2}s b_{2i-1}^{2s-1})\\
&=& -2(\gamma_{\ell} + a t) (1+c)b_{2i}^{2s}.
\end{eqnarray*}
Similarly we can prove that if the theorem holds for $t \leq 2s+1$, then it also holds for $t=2s+2$.
\end{proof}

The previous proof can be used to give explicit formulae for the coefficients $b_j^i$ in the expansion of $\psi_{t,\ell}$. This yields the following result.
\begin{theorem}[Explicit form]
The coefficients $b_{i}^{j}$ in the expansion of the functions $\psi_{t,\ell}$ take the following form
\begin{eqnarray*}
b_{2i}^{2t} &=& 2^{2t} (1+c)^{2t} \left( \begin{array}{l}t\\i \end{array} \right)\frac{\Gamma(\gamma_{\ell}/a + t)}{\Gamma(\gamma_{\ell}/a + i)}\left( \frac{a}{2}\right)^{t-i}\\
b_{2i+1}^{2t+1} &=&  -2^{2t+1} (1+c)^{2t+1}\left( \begin{array}{l}t\\i \end{array} \right)\frac{\Gamma(\gamma_{\ell}/a + t+1)}{\Gamma(\gamma_{\ell}/a + i+1)}\left( \frac{a}{2}\right)^{t-i}.
\end{eqnarray*}
\label{explformhermite}
\end{theorem}

\begin{proof}
We first prove the formula for $b_{2i}^{2t}$. Using the expressions from the previous proof we obtain
\begin{eqnarray*}
b_{2i}^{2t} &=& -2 (1+c) \frac{t}{i}b_{2i-1}^{2t-1}\\
&=& 4(1+c)^{2}\frac{t}{i}\frac{(\gamma_{\ell} +a(t-1))}{(\gamma_{\ell} +a(i-1))}b_{2i-2}^{2t-2}\\
&=& \ldots\\
&=&2^{2i}(1+c)^{2i} \frac{t\ldots (t-i+1)}{i(i-1)\ldots 1}\frac{(\gamma_{\ell} +a(t-1))\ldots (\gamma_{\ell} +a(t-i))}{(\gamma_{\ell} +a(i-1)) \ldots \gamma_{\ell}}b_{0}^{2t-2i}\\
&=&2^{2i} \left( \begin{array}{l}t\\i \end{array} \right)\frac{\Gamma(\gamma_{\ell}/a + t)\Gamma(\gamma_{\ell}/a )}{\Gamma(\gamma_{\ell}/a + t-i)\Gamma(\gamma_{\ell}/a + i)}b_{0}^{2t-2i}.
\end{eqnarray*}
So we need a formula for $b_{0}^{2t}$. This can be done as follows
\begin{eqnarray*}
b_{0}^{2t} &=& -(1+c) \gamma_{\ell} b_{1}^{2t-1}\\
&=&2(1+c)^{2} (\gamma_{\ell} + a(t-1)) b_{0}^{2t-2}\\
&=&\ldots\\
&=&2^{t}(1+c)^{2t} a^{t} \frac{\Gamma(\gamma_{\ell}/a + t)}{\Gamma(\gamma_{\ell}/a)}.
\end{eqnarray*}
Combining these results gives the desired formula for $b_{2i}^{2t}$. The formula for $b_{2i+1}^{2t+1}$ follows from the observation that
\[
b_{2i+1}^{2t+1} = -2(1+c)\frac{\gamma_{\ell} + a t}{\gamma_{\ell} + a i} b_{2i}^{2t}.
\]
\end{proof}

We can now connect the functions $\psi_{t,\ell}$ with Laguerre polynomials on the real line. This is the topic of the next theorem.

\begin{theorem}
One has that
\begin{eqnarray*}
\psi_{2t,\ell} &=& 2^{2t} (1+c)^{2t} t! \left( \frac{a}{2} \right)^{t} L_{t}^{\frac{\gamma_{\ell}}{a} -1}\left(\frac{2}{a}r^{a}\right)r^{\beta_{\ell}}M_{\ell}\\
\psi_{2t+1,\ell} &=& -2^{2t+1} (1+c)^{2t+1} t! \left( \frac{a}{2} \right)^{t} L_{t}^{\frac{\gamma_{\ell}}{a}}\left(\frac{2}{a}r^{a}\right) \ux_{a}r^{\beta_{\ell}}M_{\ell},
\end{eqnarray*}
where $L_{n}^{\alpha}$ are the generalized Laguerre polynomials on the real line.
\label{LaguerreRelations}
\end{theorem}

\begin{proof}
This follows immediately by comparing the coefficients given in theorem \ref{explformhermite} with the definition of the generalized Laguerre polynomials:
\[
L_{t}^{\alpha}(x) = \sum_{i=0}^t \frac{\Gamma(t +\alpha +1)}{i! (t-i)! \Gamma(i + \alpha +1)} (-x)^i.
\]
\end{proof}

Now we introduce the following $\cC l_{0,m}$-valued inner product
\begin{equation}
\langle f, g \rangle = \int_{\mR^{m}} \overline{f} g \;h(r) w_{k}(x) dx
\end{equation}
where $h(r)$ is the measure associated to $\dD$ (see proposition \ref{partIntProp}) and with $\bar{.}$ the main anti-involution on the Clifford algebra $\cC l_{0,m}$ defined by
\begin{eqnarray*}
\overline{a b} &=& \overline{b} \overline{a}\\
\overline{e_{i}} &=& -e_{i}, \quad i = 1,\ldots, m.
\end{eqnarray*}
The set of functions $\phi_{t,\ell}$ satisfies nice orthogonality relations with respect to this inner product. They are given in the following theorem.
\begin{theorem}[Orthogonality]
One has
\begin{equation}
\langle \phi_{t,\ell} , \phi_{s,m} \rangle = c(t,\ell) \delta_{t s} \delta_{\ell m} \int_{\mS^{m-1}} \overline{M_{\ell}}(\xi) M_{\ell}(\xi) d\sigma(\xi)
\end{equation}
where
\begin{eqnarray*}
c(2t,\ell) &=& \frac{1}{2} (2a)^{2t} (1+c)^{4t} t! \Gamma\left(\frac{\gamma_{\ell}}{a} + t \right) \left(\frac{a}{2}\right)^{\frac{\gamma_{\ell}}{a}}\\
c(2t+1,\ell) &=& \frac{1}{2} (2a)^{2t+1} (1+c)^{4t+2} t! \Gamma\left(\frac{\gamma_{\ell}}{a} + t + 1\right) \left(\frac{a}{2}\right)^{\frac{\gamma_{\ell}}{a}}.
\end{eqnarray*}
\label{orthLaguerre}
\end{theorem}

\begin{proof}
There are two possibilities to obtain this result. One can use the expression of $\phi_{t,\ell}$ and $\phi_{s,m}$ in terms of Laguerre polynomials (see theorem \ref{LaguerreRelations}) and then reduce this to the well-known orthogonality relation of the Laguerre polynomials on the real line, combined with the orthogonality of Dunkl monogenics of different degree on the unit sphere (see \cite{DBGeg}). 

Alternatively, one can note that the adjoint of $\dD - (1+c)\ux_{a}$ with respect to $\langle,\rangle$ is given by $\dD + (1+c)\ux_{a}$. 
We then have the following calculation (suppose $t \geq s$)
\begin{eqnarray*}
\langle \phi_{t,\ell} , \phi_{s,m} \rangle  &=& \langle \left(\dD - (1+c) \ux_{a} \right) \phi_{t-1,\ell} , \phi_{s,m} \rangle\\
&=& \langle  \phi_{t-1,\ell} , \left(\dD + (1+c) \ux_{a} \right)\phi_{s,m} \rangle\\
&=& \langle  \phi_{t-1,\ell} , e^{-r^{a}/a}\dD \psi_{s,m} \rangle\\
&=& C(s,m) \langle  \phi_{t-1,\ell} , e^{-r^{a}/a}\psi_{s-1,m} \rangle\\
&=& \ldots \\
&=& C(s,m) \ldots C(1,m)\langle  \phi_{t-s,\ell} , e^{-r^{a}/a}\psi_{0,m} \rangle\\
&=&C(s,m) \ldots C(1,m) \delta_{t s} \delta_{\ell m} \langle  \phi_{0,\ell} , \phi_{0,\ell} \rangle
\end{eqnarray*}
where we used theorem \ref{DiffequationCHpolys} and the orthogonality of Dunkl monogenics of different degree (see \cite{DBGeg}). 
Finally, we obtain
\begin{eqnarray*}
\langle  \phi_{0,\ell} , \phi_{0,\ell} \rangle &=& \langle  r^{\beta_{\ell}}M_{\ell} e^{-r^{a}/a} , r^{\beta_{\ell}}M_{\ell} e^{-r^{a}/a} \rangle\\
&=& \frac{1}{2} \left( \frac{a}{2}\right)^{\frac{\gamma_{\ell}}{a}}\Gamma(\gamma_{\ell}/a) \int_{\mS^{m-1}} \overline{M_{\ell}}(\xi) M_{\ell}(\xi) d\sigma(\xi).
\end{eqnarray*}
Putting everything together and substituting the values of $C(s,m)$ then yields the normalization constants.
\end{proof}

We can also associate a quantum harmonic oscillator with the operator $\dD$. This will be the basis for developing a Fourier transform in the next section. The equation of this harmonic oscillator is discussed in the following theorem.

\begin{theorem}[Harmonic oscillator]
The functions $\phi_{t,\ell}$ satisfy the following second-order PDE
\begin{equation}
\left(\dD^{2} - (1+c)^{2}\ux_{a}^{2}\right) \phi_{t,\ell} = (1+c)^{2}(\gamma_{\ell} + a t ) \phi_{t,\ell}.
\end{equation}
\label{HarmOsc}
\end{theorem}

\begin{proof}
Using the gauge property of $\dD$ (see lemma \ref{gaugeD}) we calculate consecutively
\begin{eqnarray*}
\dD^{2}\phi_{t,\ell} &=& \dD^{2} \psi_{t,\ell} e^{-r^{a}/a}\\
&=&e^{-r^{a}/a} \left(\dD - (1+c)\ux_{a} \right)^{2} \psi_{t,\ell}\\
&=& e^{-r^{a}/a} \left(\dD^{2} - 2 (1+c) \ux_{a}\dD+ (1+c)^{2}\ux_{a}^{2} - (1+c) [\dD, \ux_{a}] \right) \psi_{t,\ell}\\
&=& e^{-r^{a}/a} \left(C(t,\ell)+ (1+c)^{2}\ux_{a}^{2} - (1+c) [\dD, \ux_{a}] \right) \psi_{t,\ell},
\end{eqnarray*}
where we have also used the differential equation satisfied by $\psi_{t,\ell}$ (see theorem \ref{DiffequationCHpolys}). To simplify this further, we need to calculate the action of $[\dD, \ux_{a}]$ on $\psi_{t,\ell}$. We first prove that $[[\dD, \ux_{a}], \ux_{a}^{2}] = 0$. Indeed, using theorem \ref{ospFamily} we obtain
\begin{eqnarray*}
[[\dD, \ux_{a}], \ux_{a}^{2}] &=& \dD \ux_{a}^{3} -\ux_{a}\dD \ux_{a}^{2} - \ux_{a}^{2}\dD \ux_{a} +\ux_{a}^{3}\dD \\
&=&\dD \ux_{a}^{3}  +\ux_{a}^{3}\dD  - \ux_{a}\{\dD, \ux_{a}\}\ux_{a}\\
&=& \ux_{a}^{2}\dD \ux_{a} -a(1+c)\ux_{a}^{2} +\ux_{a}^{3}\dD +2(1+c)\ux_{a}\left(\mE + \frac{\delta}{2} \right)\ux_{a}\\
&=& -2(1+c) \ux_{a}^{2} \left(\mE + \frac{\delta}{2} \right) -a(1+c)\ux_{a}^{2}+2(1+c)\ux_{a}\left(\mE + \frac{\delta}{2} \right)\ux_{a}\\
&=&-2(1+c) \ux_{a}[\ux_{a}, \mE] -a(1+c)\ux_{a}^{2}\\
&=&0.
\end{eqnarray*}
Using this result, combined with theorem \ref{LaguerreRelations}, we find
\begin{eqnarray*}
[\dD, \ux_{a}] \psi_{2t,\ell} &=&2^{2t} (1+c)^{2t} t! \left( \frac{a}{2} \right)^{t} L_{t}^{\frac{\gamma_{\ell}}{a} -1}\left(\frac{2}{a}r^{a}\right) [\dD, \ux_{a}] r^{\beta_{\ell}}M_{\ell}\\
&=&- (1+c) \gamma_{\ell} \, 2^{2t} (1+c)^{2t} t! \left( \frac{a}{2} \right)^{t} L_{t}^{\frac{\gamma_{\ell}}{a} -1}\left(\frac{2}{a}r^{a}\right) r^{\beta_{\ell}}M_{\ell}\\
&=&- (1+c) \gamma_{\ell} \,  \psi_{2t,\ell}
\end{eqnarray*}
and similarly
\[
[\dD, \ux_{a}] \psi_{2t+1,\ell} = -(1+c) (a-\gamma_{\ell}) \psi_{2t+1,\ell}.
\]
Combining these results with the previous calculation then yields the theorem.
\end{proof}

\section{Associated Fourier transforms}
\label{SecFourier}

Denote by $M_{\ell}^{(m)} \in \cM_{\ell}$, ($m = 1,\ldots, \dim \cM_{\ell}$), an orthonormal basis of the space of Dunkl monogenics of degree $\ell$, in the sense that
\[
\int_{\mS^{m-1}} \left[\overline{M_{\ell}}^{(m_{1})}(\xi) M_{\ell}^{(m_{1})} (\xi) \right]_{0} d\sigma(\xi) = \delta_{m_{1} m_{2}}
\]
where $[\, . \, ]_{0}$ denotes the projection on the 0-vector part in the Clifford algebra $\cC l_{0,m}$. Then the functions 
\[
\phi_{t, \ell, m} = (c(t, \ell))^{-1/2} \left(\left(\dD - 2(1+c) \ux_{a} \right)^{t} r^{\beta_{\ell}} M_{\ell}^{(m)}\right) e^{-r^{a}/a} 
\]
with $t, \ell \in \mN$ and $m = 1,\ldots, \dim \cM_{\ell}$ form an orthonormal basis for the Hilbert space $L_{2}(\mR^{m}, h(r) w_{k}(x) dx)$ of functions taking values in $\cC l_{0,m}$, equipped with the inner product
\[
\langle f, g \rangle = \int_{\mR^{m}} \left[\, \overline{f} g \, \right]_{0} \;h(r) w_{k}(x) dx.
\]

Note that theorem \ref{HarmOsc} can now be rewritten as follows 
\[
\frac{1}{ a(1+c)^{2}}\left(\dD^{2} - (1+c)^{2}\ux_{a}^{2}\right) \phi_{t,\ell,m} - \left(\frac{1}{2}+ \frac{\mu - 1}{a(1+c)}\right) \phi_{t,\ell,m} =(\frac{2\ell}{a(1+c)} + t ) \phi_{t,\ell,m}.
\]
We can hence introduce the associated Fourier transform by
\begin{equation}
\cF_{\dD} = e^{i \frac{\pi}{2} \left(\frac{1}{2}+ \frac{\mu - 1}{a(1+c)} \right)} e^{\frac{-i \pi}{ 2a(1+c)^{2}}\left(\dD^{2} - (1+c)^{2}\ux_{a}^{2}\right)}.
\end{equation}
It is clear that this transform satisfies
\[
\cF_{\dD}(\phi_{t,\ell,m}) = (-i)^{t} e^{-\frac{i \pi \ell}{a(1+c)}} \phi_{t,\ell,m}.
\]
Using the Hadamard identity for two linear maps $A,B$ given by
\[
\exp(A)B\exp(-A)=ad_{A}B=B+[A,B]+\frac{1}{2!}[A,[A,B]]+\dots,
\]
one can check that
\begin{eqnarray*}
\cF_{\dD}(\dD \, \cdot) &=& i (1+c) \ux_{a} \cF_{\dD} (\cdot)\\
\cF_{\dD}(\ux_{a} \, \cdot) &=& \frac{i}{1+c} \dD \cF_{\dD}(\cdot).
\end{eqnarray*}
If we represent the Fourier transform $\cF_{\dD}$ as an integral transform as follows
\[
\cF_{\dD}(f) = \int_{\mR^{m}} K(x,y) \; f(x) \; h(r) w_{k}(x) dx
\]
with $K(x,y)$ the integral kernel, then we find that, using proposition \ref{partIntProp}, the kernel has to satisfy the system of PDEs given by
\begin{eqnarray*}
i \dD_{y} K &=& (1+c) K\ux_{a}\\
i K \dD_{x} &=& (1+c) \uy_{a} K,
\end{eqnarray*}
where the subscript denotes the variables under consideration.


Note that, as in general $\dD^{2}$ is not a scalar operator, we do not expect the kernel $K(x,y)$ of the associated Fourier transform to be a scalar function, but rather a $\cC l_{0,m}$-valued function.

\subsection{Determination of the kernel when $\dD^{2}$ is scalar}

Let us now consider the special case where $\dD^{2}$ is scalar. Then $c = \frac{2}{a}-1$ and $\dD$ reduces to
\[
\dD = r^{1- \frac{a}{2}}\cD_{k} + b r^{-\frac{a}{2}-1} \ux + \left(\frac{2}{a}-1 \right)r^{-\frac{a}{2}-1}\ux \mE.
\]
The Fourier transform now acts as follows
\[
\cF_{\dD}(\phi_{t,\ell,m}) = (-i)^{t + \ell}\phi_{t,\ell,m}
\]
which is the same behaviour as the classical Fourier or Dunkl transform.

As the Fourier transform is proportional to
\[
e^{\frac{-i \pi a}{ 8}\left(\dD^{2} - \frac{4}{a^{2}}\ux_{a}^{2}\right)},
\]
we also know that $K(x,y)$ is a scalar function and that this function is determined by the system ($j = 1, \ldots,m$)
\begin{equation}
\left(r_{x}^{1-\frac{a}{2}} T_{j,x} + b r_{x}^{-1-\frac{a}{2}} x_{j} + \left(\frac{2}{a}-1 \right) r_{x}^{-\frac{a}{2}} x_{j} \partial_{r_{x}}\right) K(x,y) = -\frac{2}{a}i y_{j}r_{y}^{\frac{a}{2}-1}K(x,y),
\label{ScalarSystemDunkl}
\end{equation}
where we again use subscripts $x$ and $y$ to denote the variables under consideration. For general reflection groups this system seems quite complicated to solve. So we start with discussing the case where the symmetry is $O(m)$ (= non-Dunkl case with $k=0$). Then the system reduces to  ($j = 1, \ldots,m$)
\begin{equation}
\left(r_{x}^{1-\frac{a}{2}} \partial_{x_{j}} + b r_{x}^{-1-\frac{a}{2}} x_{j} + \left(\frac{2}{a}-1 \right) r_{x}^{-\frac{a}{2}} x_{j} \partial_{r_{x}}\right) K(x,y) = -\frac{2}{a} i y_{j}r_{y}^{\frac{a}{2}-1}K(x,y).
\label{ScalarSystem}
\end{equation}
Multiplying each equation by $x_{j}$ and summing from $1$ to $m$ yields
\[
r_{x}\partial_{r_{x}} K + \frac{ab}{2} K = - i \langle \ux, \uy \rangle (r_{x} r_{y})^{ \frac{a}{2}-1} K.
\]
This equation can be integrated immediately, leading to
\begin{equation}
K(x,y) = d (r_{x} r_{y})^{ -\frac{ab}{2}} e^{-\frac{2i}{a}\langle \ux, \uy \rangle (r_{x} r_{y})^{ \frac{a}{2}-1}}
\label{exactKernelFourier}
\end{equation}
where $d$ is a constant that is still to be determined. Note that it can be checked that (\ref{exactKernelFourier}) indeed satisfies the system (\ref{ScalarSystem}). To determine the constant, we put $f = \phi_{0,0} = r_{x}^{-\frac{ab}{2}} e^{-r_{x}^{a}/a}$ and calculate
\[
\int_{\mR^{m}} K(x,y) f(x) h(r) dx
\]
where we have put $w_{k}(x) =1$ as we are considering the non-Dunkl case. Using the Funk-Hecke theorem and an integral identity for Bessel functions, we can calculate this integral as
\[
\int_{\mR^{m}} K(x,y) f(x) h(r) dx = d (2\pi)^{\frac{m}{2}} \left( \frac{2}{a} \right)^{1-\frac{m}{2}} r_{y}^{-\frac{ab}{2}} e^{-r_{y}^{a}/a}
\]
so we conclude that
\[
d = (2\pi)^{-\frac{m}{2}} \left( \frac{2}{a}\right)^{\frac{m}{2}-1}.
\]

Hence, in the case $k=0$ we have obtained that
\begin{equation}
\cF_{\dD} (g) = (2\pi)^{-\frac{m}{2}} \left( \frac{2}{a}\right)^{\frac{m}{2}-1} \int_{\mR^{m}}(r_{x} r_{y})^{ -\frac{ab}{2}} e^{-\frac{2i}{a}\langle \ux, \uy \rangle (r_{x} r_{y})^{ \frac{a}{2}-1}} g(x) h(r_{x}) dx
\label{ExplicitFourier}
\end{equation}
with the measure $h(r_{x}) = r_{x}^{ab - (1-\frac{a}{2})m}$.

\vspace{5mm}
\textbf{Interpretation:} 

By slightly changing the Laplace operator proposed in \cite{Orsted2} from $r^{2-a}\Delta_{k}$ to
\begin{eqnarray*}
\dD^{2} &=&-r^{2-a} \Delta_{k} -b\left( b+  \mu -2 \right) r^{-a} - \left(4 \frac{b}{a} + \left(\frac{2}{a}-1\right)^{2} + \left(\frac{2}{a}-1\right)\mu \right) r^{1-a}\partial_{r}\\
&& -  \left(\frac{4}{a^{2}}-1 \right) r^{2-a} \partial_{r}^{2},
\end{eqnarray*}
we have obtained a Fourier transform which behaves exactly as the classical Fourier transform (i.e. that has the same eigenvalues) and of which the kernel is known explicitly. Moreover, this modified Laplace operator can also be factorized ($\dD^{2} = - \sum_{i=1}^{m} D_{i}^{2}$). Note that in the resulting operator there is still quite a lot of freedom, as there are no restrictions on the value of $b$. So if one takes e.g. $b = 0$, the result simplifies even more.

\vspace{5mm}
To gain more insight in this new Fourier transform, we introduce the following two operators $P$ and $Q$ defined by
\begin{eqnarray*}
P f(x) &=& r^{b} f\left( \left(\frac{a}{2} \right)^{\frac{1}{a}}\ux r^{\frac{2}{a}-1}\right)\\
Q f(x) &=& r^{-\frac{ab}{2}} f\left( \left(\frac{2}{a} \right)^{\frac{1}{2}}\ux r^{\frac{a}{2}-1}\right).
\end{eqnarray*}
These two operators act as generalized Kelvin transformations (for a discussion of the classical Kelvin transformation in harmonic analysis we refer the reader to \cite{MR1805196}). Indeed, it is easily calculated that
\[
Q P = P Q =  \left(\frac{2}{a} \right)^{\frac{b}{2}},
\]
hence they are each others inverse up to a constant. They allow to write the transform (\ref{ExplicitFourier}) in terms of the classical Fourier transform as follows (see also the subsequent theorem \ref{ExpFourierScalarD})
\[
\cF_{\dD} = \left(\frac{a}{2} \right)^{\frac{b}{2}}Q_{y} \cF_{\mbox{class}} P_{x},
\]
where the subscript denotes on which variables the operators act.

First we show that the Dunkl operators $T_{i}$ transform nicely under the action of $P$ and $Q$. This is the subject of the following proposition.
\begin{proposition}[Intertwining relations]
One has that
\[ 
r^{1- \frac{a}{2}}T_{i} + b r^{-\frac{a}{2}-1} x_{i} + \left(\frac{2}{a}-1 \right)r^{-\frac{a}{2}}x_{i}\partial_{r} = \left(\frac{a}{2} \right)^{\frac{b-1}{2}} Q T_{i} P, \qquad i = 1, \ldots, m. 
\]
\label{KelvinIntertwining}
\end{proposition}

\begin{proof}
We first calculate $\partial_{x_{i}} P$ as follows
\begin{eqnarray*}
&&\partial_{x_{i}} P f(x)\\ &=& T_{i}r^{b} f\left( \left(\frac{a}{2} \right)^{\frac{1}{a}}\ux r^{\frac{2}{a}-1}\right)\\
&=& br^{b-2} x_{i}  f\left( \left(\frac{a}{2} \right)^{\frac{1}{a}}\ux r^{\frac{2}{a}-1}\right) + r^{b} \partial_{x_{i}} f\left( \left(\frac{a}{2} \right)^{\frac{1}{a}}\ux r^{\frac{2}{a}-1}\right)\\
&=&br^{b-2} x_{i}  f\left( \left(\frac{a}{2} \right)^{\frac{1}{a}}\ux r^{\frac{2}{a}-1}\right) \\&&+ r^{b} \sum_{j=1}^{m}(\partial_{x_{j}}f)\left( \left(\frac{a}{2} \right)^{\frac{1}{a}}\ux r^{\frac{2}{a}-1}\right) \left(\delta_{ij} \left(\frac{a}{2} \right)^{\frac{1}{a}} r^{\frac{2}{a}-1}  + (\frac{2}{a}-1)  \left(\frac{a}{2} \right)^{\frac{1}{a}} r^{\frac{2}{a}-3}x_{i}x_{j}\right)\\
&=&br^{b-2} x_{i}  f\left( \left(\frac{a}{2} \right)^{\frac{1}{a}}\ux r^{\frac{2}{a}-1}\right)  +  \left(\frac{a}{2} \right)^{\frac{1}{a}} r^{b+\frac{2}{a}-1} (\partial_{x_{i}}f)\left( \left(\frac{a}{2} \right)^{\frac{1}{a}}\ux r^{\frac{2}{a}-1}\right)\\
&& + (\frac{2}{a}-1)  \left(\frac{a}{2} \right)^{\frac{1}{a}} r^{b+\frac{2}{a}-3}x_{i} \sum_{j=1}^{m}x_{j}(\partial_{x_{j}}f)\left( \left(\frac{a}{2} \right)^{\frac{1}{a}}\ux r^{\frac{2}{a}-1}\right).
\end{eqnarray*}
Applying $Q$ then yields
\[
\left(\frac{a}{2} \right)^{\frac{b-1}{2}}Q\partial_{x_{i}} P f(x) = \left( r^{1- \frac{a}{2}}\partial_{x_{i}} + b r^{-\frac{a}{2}-1} x_{i} + \left(\frac{2}{a}-1 \right)r^{-\frac{a}{2}}x_{i}\partial_{r} \right) f(x). 
\]
We also have
\begin{eqnarray*}
&&(T_{i }- \partial_{x_{i}} )P f(x)\\&=&(T_{i }- \partial_{x_{i}} )r^{b} f\left( \left(\frac{a}{2} \right)^{\frac{1}{a}}\ux r^{\frac{2}{a}-1}\right)\\
&=& r^{b} \sum_{\alpha \in R_{+}} k_{\alpha} \alpha_{i} \frac{f\left( \left(\frac{a}{2} \right)^{\frac{1}{a}}\ux r^{\frac{2}{a}-1}\right) - f\left( \left(\frac{a}{2} \right)^{\frac{1}{a}}r_{\alpha}(\ux) r^{\frac{2}{a}-1}\right)}{\langle \alpha,  \left(\frac{a}{2} \right)^{\frac{1}{a}}\ux r^{\frac{2}{a}-1} \rangle} \left( \left(\frac{a}{2} \right)^{\frac{1}{a}}r^{\frac{2}{a}-1}\right).
\end{eqnarray*}
Again applying $Q$ gives
\[
\left(\frac{a}{2} \right)^{\frac{b-1}{2}}Q(T_{i }- \partial_{x_{i}} )P f(x) = r^{1- \frac{a}{2}} (T_{i}-\partial_{x_{i}}) f(x)
\]
and putting everything together completes the proof.
\end{proof}

\begin{remark}
As a consequence of proposition \ref{KelvinIntertwining} we obtain that $\dD = \left(\frac{a}{2} \right)^{\frac{b-1}{2}} Q \cD_{k} P$. This means that $P$ and its inverse $Q$ act as analogues of the intertwining operator $V_{k}$ in the theory of Dunkl operators (see e.g. \cite{MR1273532}).
\end{remark}

Using proposition \ref{KelvinIntertwining} we can also give a more explicite expression for the related Fourier transform in the Dunkl case (i.e. $k \neq 0$). This the subject of the following theorem.

\begin{theorem}
When $c = 2/a-1$ and $k \geq 0$, one has the following explicit expressions of the Fourier transform related to $\dD$
\begin{eqnarray*}
\cF_{\dD} &=& e^{\frac{i \pi \mu}{4}} e^{\frac{-i \pi a}{ 8}\left(\dD^{2} - \frac{4}{a^{2}}\ux_{a}^{2}\right)}\\
&=& \left(\frac{a}{2} \right)^{\frac{b}{2}} Q_{y} \cF_{k} P_{x}(.)\\ 
&=& c_{k}^{-1} \left(\frac{2}{a} \right)^{\frac{\mu}{2}-1}  \int_{\mR^{m}} (r_{x}r_{y})^{-\frac{ab}{2}}D\left( \left(\frac{2}{a} \right)^{\frac{1}{2}}\ux r_{x}^{\frac{a}{2}-1}, -\left(\frac{2}{a} \right)^{\frac{1}{2}} i \uy r_{y}^{\frac{a}{2}-1} \right) \; (.)\;  \\&&\times h(r_{x}) w_{k}(x) dx
\end{eqnarray*}
where $D(.,.)$ is the Dunkl kernel and $h(r_{x}) = r_{x}^{ab - (1-\frac{a}{2})\mu}$.
\label{ExpFourierScalarD}
\end{theorem}

\begin{proof}
Consider the basis of $L_{2}(\mR^{m}, h(r) w_{k}(x) dx)$ given by $\{\phi_{t,\ell,m} \}$ ($t, \ell \in \mN$ and $m = 1,\ldots, \dim \cM_{\ell}$). Using the explicit expressions in theorem \ref{LaguerreRelations}, it is easy to see that the set $\{P_{x}\phi_{t,\ell,m}\}$ is an eigenfunction basis of the Dunkl transform $\cF_{k}$, i.e.
\[
\cF_{k}(P_{x}\phi_{t,\ell,m}(x)) = (-i)^{t+\ell}P_{y}\phi_{t,\ell,m}(y).
\]
Acting with $Q_{y}$ then yields the first equality in the theorem.

Now we can calculate
\begin{eqnarray*}
 Q_{y} \cF_{k} P_{x}f(x) &=& c_{k}^{-1}Q_{y} \int_{\mR^m} \,D(x, -iy)\, r_{x}^{b} f\left( \left(\frac{a}{2} \right)^{\frac{1}{a}}\ux r_{x}^{\frac{2}{a}-1}\right) w_k(x) dx\\
 &=& c_{k}^{-1} \int_{\mR^m}r_{y}^{-\frac{ab}{2}} \,D\left(x, -i\left(\frac{2}{a} \right)^{\frac{1}{2}}\uy r_{y}^{\frac{a}{2}-1}\right)\, r_{x}^{b}\\
 && \times f\left( \left(\frac{a}{2} \right)^{\frac{1}{a}}\ux r_{x}^{\frac{2}{a}-1}\right) w_k(x) dx.
\end{eqnarray*}
We perform a change of variables given by $z_{i} =  \left(\frac{a}{2} \right)^{\frac{1}{a}}x_{i} r_{x}^{\frac{2}{a}-1}$, $i=1, \ldots, m$. The absolute value of the determinant of the Jacobian of this transformation is given by 
\[
|\det{J}| =  \left(\frac{2}{a}\right)^{\frac{m  }{2}-1}r_{z}^{m(\frac{a}{2}-1)},
\]
as can be calculated using the matrix determinant lemma. Putting everything together then yields the second equality of the theorem.
\end{proof}

\subsection{The case $a = -2$}
\label{KelvinSection}

In this case, the Dirac operator takes the form
\[
\cD_{k,-2} =  \sum_{i=1}^{m} D_{i}e_{i} = r^{2}\cD_{k} - \left(\mu-2 \right) \ux -2 \ux \mE
\]
Its components $D_{i}$ are commutative and given by
\[
D_{i} = r^{2} T_{i} + (2-\mu) x_{i} - 2 r x_{i} \partial_{r}.
\]
We also introduce the operator $\textbf{I}_{k}$, acting on functions defined in $\mR^{m}$ as
\begin{equation}
\textbf{I}_{k} \left( f(x) \right) = r^{2-\mu} f\left( \frac{x}{r^{2}}\right).
\label{Kelvin}
\end{equation}
Note that if $k=0$ this operator reduces to the classical Kelvin inversion. It is moreover easy to calculate that $\textbf{I}_{k}^{2}= \mbox{id}$. Now we obtain the following proposition (the proof is similar to the proof of proposition \ref{KelvinIntertwining}):

\begin{proposition}
The operator $\textbf{I}_{k}$ intertwines the Dirac operator for $a=2$ with the Dirac operator for $a=-2$ in the following way:
\[
\cD_{k,-2} = \textbf{I}_{k}\cD_{k} \textbf{I}_{k}.
\]
\end{proposition}

We can also connect a Fourier transform with the Dirac operator $\cD_{k,-2}$. This transform is given by
\begin{eqnarray*}
\cF_{k,-2}(f) &=& \textbf{I}_{k,y} \cF_{k} \textbf{I}_{k,x} f\\
&=&c_{k}^{-1}\int_{\mR^{m}} \left(r_{x}r_{y}\right)^{2-\mu}  D\left( \frac{x}{r_{x}^{2}}, -i\frac{y}{r_{y}^{2}} \right) \; f(x) \; r_{x}^{-4} w_{k}(x)dx.
\end{eqnarray*}

\section{Conclusions and outlook for further research}
\label{SecConclusions}

In this paper we have obtained a new class of deformations of the realization of $\mathfrak{osp}(1|2)$ using Dunkl operators. Some elements of this new class clearly have preferred properties. In the case where $c = 2/a-1$ the square of $\dD$ is scalar and one can explicitly construct the related Fourier transform and the intertwining operator. 
Another interesting case appears to be when $b=c=0$. This is in some sense the easiest deformation of the classical $\mathfrak{osp}(1|2)$ realization and in spirit the closest to the $\mathfrak{sl}_{2}$ deformation obtained in \cite{Orsted2}. One expects that a better understanding of this case will also yield the key to study the more general cases.
Finally, also the deformation in formula (\ref{aDirac}) seems interesting, as this case preserves as much as possible of the $\mathfrak{sl}_{2}$ relations already obtained in \cite{Orsted2}.

It is important to note that the class of deformations we have obtained is in some sense relatively limited. It is in general possible to add additional terms to $\dD$, as long as they are of homogeneity $-a$. Typical examples of such terms would be powers of $\ux \cD_{k} \mE$ multiplied with suitable powers of $r$. This would allow to introduce even more parameters, and at this point it is not entirely clear whether one can still study the resulting operators in a systematic way.

For some special values of the parameters, we were able to construct an explicit intertwining map relating the deformed Dirac operator with the classical one. One can not expect to obtain explicit constructions of intertwinors between all elements of the family (as e.g. in the theory of Dunkl differential operators the intertwining operators are only explicitly known for a few special cases). However, the existence or non-existence of such operators is clearly a very important topic for further research.

In the case of the deformation of the realization of $\mathfrak{sl}_{2}$ introduced in \cite{Orsted2}, the parameter $a$ interpolates between representations of different Lie algebras. From that point of view, one expects something similar to happen for the family of operators introduced in this paper. Our expectation is that one of them is realized in the kernel of Dirac operator (which is similar to the realization of the minimal orthogonal representation in the kernel of the conformally invariant Laplace-Yamabe operator). It is subject of ongoing work to describe the kernel of the conformally invariant Dirac operator (see \cite{Orsted3}).
More abstractly, in \cite{Orsted2} everything is based on the emergence of dual pairs with one member $\mathfrak{sl}_{2}$ and the other some Lie algebra. Similarly, in our case one should consider a theory of super dual pairs in Lie superalgebras, one of the members being $\mathfrak{osp}(1|2)$.

It would be very insightful if one could obtain explicit expressions of the kernel of the Fourier transform related to $\dD$ for other special values of the parameters. At least in the non-Dunkl case ($k =0$) one hopes that the first order system describing the kernel can somehow be solved explicitly. This ties in with the general study of Fourier transforms in Clifford analysis inaugurated in \cite{DBFourier} and might also shed new light on the Fourier transforms introduced in \cite{Orsted2}.
A related question is whether there exists an uncertainty relation for the new type of Fourier transforms.
Also, so far we have only studied the special case $z=- i \pi/2$ of the holomorphic semigroup
\[
e^{\frac{z}{ a(1+c)^{2}}\left(\dD^{2} - (1+c)^{2}\ux_{a}^{2}\right)}.
\]
In the future we plan to study this semigroup for general complex values $z$.

\section*{Acknowledgement}
The authors would like to thank Yuan Xu for helpful discussions regarding this preprint.

\end{document}